\documentclass[a4paper]{article}

\usepackage[latin1]{inputenc}
\usepackage[T1]{fontenc}
\usepackage[francais,english]{babel}

\usepackage{amsmath,amssymb,amsfonts,amsthm}
\usepackage{amsfonts}
\usepackage{amssymb}
\usepackage{latexsym}
\usepackage{comment}
\usepackage{color}
\usepackage{eqsection}

\def\R{\mathbb R}

\def\Z{\mathbb Z}
\def\C{\mathbb C}

\def\I{\mathbb I}

\def\T{\mathbb T}

\def\et0{e^{tA}x_0}

\def\hl{A}
\def\sph{H}
\def\spv{V}

\newtheorem{thm}{Theorem}
\newtheorem{Prop}{Proposition}
\newtheorem{Def}{Definition}
\newtheorem{Lem}{Lemma}

\newtheorem{rk}{Remark}

\author{
K. \textsc{Beauchard} 
\footnote{ENS Rennes, Avenue Schumann, 35170 Bruz, France,
email: Karine.Beauchard@ens-rennes.fr (corresponding author)},
P. \textsc{Cannarsa}
\footnote{Dipartimento di Matematica, Universit\`a di Roma Tor Vergata, via della Ricerca Scientifica 1, 00133, Roma, Italy, 
email: cannarsa@axp.mat.uniroma2.it}
\thanks{This research has been performed in the framework of the GDRE CONEDP.  The first author was partially supported by 
the ``Agence Nationale de la Recherche'' (ANR) Projet Blanc EMAQS number ANR-2011-BS01-017-01.}
}

\title{Heat equation on the Heisenberg group: observability and applications.}
\date{}

\begin{document}

\maketitle

\begin{abstract}
We investigate observability and Lipschitz stability 
for the Heisenberg heat equation on the rectangular domain 
$$\Omega = (-1,1)\times\mathbb{T}\times\mathbb{T}$$
taking as observation regions  slices of the form $\omega=(a,b) \times \mathbb{T} \times \mathbb{T}$ 
or tubes $\omega = (a,b) \times \omega_y \times \mathbb{T}$, with $-1<a<b<1$.
We prove that observability fails for an arbitrary time $T>0$ but 
both observability and Lipschitz stability hold true after a positive minimal time, which depends on the distance between $\omega$ and the boundary of $\Omega$:
$$T_{\min} \geqslant \frac{1}{8} \min\{(1+a)^2,(1-b)^2\}.$$
Our proof follows a mixed strategy which combines the approach by Lebeau and Robbiano, 
which relies on Fourier decomposition, with Carleman inequalities for the heat equations that are solved by the Fourier modes. We  extend the  analysis to the unbounded domain $(-1,1)\times\mathbb{T}\times\mathbb{R}$.

\end{abstract}



\section{Introduction}

This article focuses on the heat equation on the Heisenberg group
\begin{equation} \label{H}
\left\lbrace \begin{array}{ll}
\big(\partial_t  - \partial_x^2 - (x \partial_z+\partial_y)^2 \big) g = \widetilde{h} & \mbox{in}\; (0,T) \times \Omega\,, \\
g(t,\pm 1,y,z)=0\,,                                                                                        & (t,y,z) \in (0,T)\times \mathbb{T} \times \mathbb{T}\,, \\
g(0,x,y,z)=g^0(x,y,z)\,,                                                                                   & (x,y,z) \in \Omega\,,
\end{array} \right.
\end{equation}
where  $\mathbb{T}$ is the $1$D-torus and $\Omega=(-1,1) \times \mathbb{T} \times \mathbb{T}$. In section~\ref{se:wp}, we will give the precise notion of weak solution to problem \eqref{H} for 
$$g^0\in L^2(\Omega)\quad\mbox{and}\quad\widetilde{h}\in L^2\big((0,T)\times\Omega\big).$$
For the above problem, we will investigate observability and Lipschitz stability. 
We recall the definition of these two notions below and we state our main results.
 
\subsection{Observability and null controllability}

\begin{Def} 
[Observability] Let $T>0$ and $\omega$ be an open subset of $\Omega$. 
System \eqref{H} is {\em observable in $\omega$ in time} $T$ if there exists a constant $C_T>0$ such that,
for every $g^0 \in L^{2}(\Omega)$, the solution of \eqref{H} with $\widetilde{h} \equiv 0$
satisfies
\begin{equation}\label{eq:observability}
\int\limits_{\Omega} |g(T,x,y,z)|^{2} dx dy dz \leqslant C_T \int\limits_{0}^{T} \int\limits_{\omega} |g(t,x,y,z)|^{2} dx dy dz dt\, .
\end{equation}
\end{Def}

\begin{thm} \label{Thm:Obs}
Let  
\begin{equation}
\label{eq:omega}
\begin{array}{c}
\omega:=(a,b) \times \omega_y \times \mathbb{T}\,,  \\
\mbox{where}  -1 < a < b < 1\; \mbox{and}\; \omega_y\; \mbox{is an open subset of} \;\mathbb{T}.
\end{array}
\end{equation} 
Then, there exists $T_{\min} \geqslant \frac{1}{8} \max\{ (1+a)^2 , (1-b)^2 \}$ such that
\begin{itemize}
\item for every $T>T_{\min}$, system \eqref{H} is observable in $\omega$ in time $T$,
\item for every $T<T_{\min}$, system \eqref{H} is not observable in $\omega$ in time $T$.
\end{itemize}
\end{thm}

It is well-known that the Heisenberg laplacian
\begin{equation}\label{eq:hl}
\hl:=-\partial_x^2 - (x \partial_z + \partial_y)^2
\end{equation}
is an hypoelliptic operator of the form $X_1^2+X_2^2$, where 
$$X_1(x,y,z):=
\left( \begin{array}{c}
1 \\ 0 \\ 0
\end{array} \right)\,, \quad \quad 
X_2(x,y,z):=
\left( \begin{array}{c}
0 \\ 1 \\ x
\end{array} \right)\,,$$
see \cite{MR0222474}. However, no clear connection between hypoellipticity and observability has been established so far.
\\

We observe that,  given the width $\ell=b-a\in (0,2)$, there is no location of the slice $\omega=(a,b)\times\mathbb{T}\times\mathbb{T}$ 
for which the minimal observability time $T_{\min}$ vanishes. 
Such a behavior differs from the one observed  for the 2D Grushin case 
$$\left\lbrace \begin{array}{ll}
\big(\partial_t  - \partial_x^2 - x^2 \partial_y^2 \big) g= 0 & \mbox{in} \;\Omega_G:= (0,T) \times (-1,1)\times(0,1)\,, \\
g(t,x,y)=0\,,                                                             & (t,y,z) \in (0,T)\times\partial \Omega_G
\\
g(0,x,y)=g^0(x,y)\,,                                                                                   & (x,y) \in \Omega_G\,,
\end{array} \right.$$
for which 
\begin{itemize}
\item $T_{\min}>0$ when $\omega=(a,b)\times(0,1)$ and $a>0$ (see \cite{Grushin}), 
\item $T_{\min}=0$  when $\omega=(0,b)\times (0,1)$ (see \cite{minimal}).
\end{itemize}
This difference may be related to the fact that, for the Heisenberg operator, the number of iterated Lie brackets of the vector fields required to generate $\mathbb{R}^3$ has no jump at $\{x=0\}$: $X_1$, $X_2$ and $[X_1,X_2]$ are needed everywhere.
\\

As usual, by the Hilbert uniqueness method (see \cite{JLL_book}, \cite{coron_book}), the observability result of Theorem \ref{Thm:Obs}
is equivalent to the following null controllability result.
\begin{Def} 
[Null controllability] Let $T>0$ and $\omega$ be an open subset of $\Omega$. 
System \eqref{H} is said to be {\em null controllable from $\omega$ in time} $T$ if,
for every $g^0 \in L^{2}(\Omega)$, there exists $\widetilde{h} \in L^2((0,T)\times\Omega)$, supported on $[0,T]\times\omega$,
such that the solution of \eqref{H} satisfies $g(T,\cdot)=0$.
\end{Def}
\begin{thm} 
Let $\omega$ be as in \eqref{eq:omega}. Then, there exists 
$$T_{\min} \geqslant \frac{1}{8} \max\{ (1+a)^2 , (1-b)^2 \}$$ such that
\begin{itemize}
\item for every $T>T_{\min}$, system \eqref{H} is null controllable from $\omega$ in time $T$,
\item for every $T<T_{\min}$, system \eqref{H} is not null controllable from $\omega$ in time $T$.
\end{itemize}
\end{thm}

\subsection{Lipschitz stability}

Taking a source term of the form
\begin{equation} \label{form_h}
\begin{array}{c}
\widetilde{h}(t,x,y,z)=R(t,x) h(x,y,z) 
\vspace{.1cm}
\\
\text{ where } R \in \mathcal C\big([0,T]\times[-1,1]\big) \text{ and } h \in L^2(\Omega),
\end{array}
\end{equation}
we will obtain Lipschitz stability estimates for \eqref{H} in the following sense.

\begin{Def}[Lipschitz stability]
Let $T>0$, let $0 \leqslant T_0 < T_1 \leqslant T$, and let $\omega$ be an open subset of $\Omega$.
We say that system \eqref{H}, with $\widetilde{h}$ as in \eqref{form_h}, satisfies a \emph{Lipschitz stability estimate on $(T_0,T_1)\times\omega$} 
if there exists a constant $\widetilde{C}_T>0$ such that, for every $g^0 \in L^2(\Omega)$ and $h \in L^2(\Omega)$,
the solution of \eqref{H} satisfies
\begin{eqnarray*}
\lefteqn{\int\limits_{\Omega} |h(x,y,z)|^2 dx dy dz }
\\
&\leqslant&
\widetilde{C}_T \Big( \int\limits_{T_0}^{T_1} \int\limits_{\omega} |\partial_t g(t,x,y,z)|^2 dx dy dz dt 
+ \int\limits_{\Omega} \left| \hl  g(T_1,x,y,z)\right|^2 dx dy dz \Big)\,,
\end{eqnarray*}
where $\hl$ is defined in \eqref{eq:hl}.
\end{Def}
Notice that the above Lipschitz stability estimate implies the uniqueness of the source term $h$
via 2 measurements: $\partial_t g|_{(T_0,T_1)\times\omega}$ and $Ag(T_1,\cdot)$.
\\

When $\omega$ is a slice, parallel to the $(y,z)$-plane, we can prove Lipschitz stability in large time under general assumptions on $R$.

\begin{thm} \label{Thm:Lips_Stab_simple}
Let $a,b \in \mathbb{R}$ be such that $-1 < a < b < 1$ and $\omega:=(a,b) \times \mathbb{T} \times \mathbb{T}$.
Suppose further that 
\begin{equation} \label{Hyp:R}
\begin{array}{c}
R, \partial_t R \in \mathcal C\big([0,T]\times[-1,1]\big) \text{ and } 
\vspace{.1cm}\\
\text{$\exists\;T_1 \in (0,T]$ and $\rho_0>0$ such that $R(T_1,x) \geqslant \rho_0, \forall x \in [-1,1]$}.
\end{array}
\end{equation}
Then, there exists $T^*>0$ such that system \eqref{H} satisfies a Lipschitz stability estimate on $(T_0,T_1)\times\omega$
 for every $T_0,T_1 \in [0,T]$ with $(T_1-T_0)>T^*.$
\end{thm}

More generally, when $\omega$ is a tube along the $z$-axis, we can still prove Lipschitz stability in large time
under an additional smallness assumption on the source term, which is probably due just to technical reasons.

\begin{thm} \label{Thm:Lips_Stab}
Let $\omega$ be as in \eqref{eq:omega}.
There exists $T^*>0$ and a continuous function $\eta:(T_*,\infty) \rightarrow (0,\infty)$ such that,
if $R$ satisfies \eqref{Hyp:R}, $T_0, T_1 \in [0,T]$, $(T_1-T_0)>T^*$, and 
\begin{equation} \label{Hyp_R_small_var}
\frac{1}{\rho_0} \Big( \int_{T_0}^{T_1} \| \partial_t R(t)\|_{L^\infty(-1,1)}^2 dt\Big)^{1\over2} < \eta(T_1-T_0),
\end{equation}
then system \eqref{H}  satisfies a Lipschitz stability estimate on $(T_0,T_1)\times\omega$.
\end{thm}

\subsection{Motivations and bibliographical comments}

\subsubsection{Motivations}

The relevance of the Heisenberg group to quantum mechanics has long been acknowledged. Indeed,  it was recognized by Weyl~\cite{Weyl}  that the Heisenberg algebra generated by the momentum and position operators comes from a Lie algebra representation associated with a corresponding group---namely the Heisenberg group (Weyl group in the traditional language of physicists).  In such a group, the role played by the so-called Heisenberg laplacian is absolutely central, being analogous to the standard laplacian in Euclidean spaces, see \cite{MR983366}. On an even larger scale, deep connections have been pointed out between the properties of subriemannian operators, like the Heisenberg laplacian, and other topics of interest to current  mathematical research such as  isoperimetric problems and systems theory, see, for instance, \cite{capogna}.

\subsubsection{Observability} 

Observability is well known to hold  for  the linear heat equation in arbitrary positive time $T$
  with any observation domain $\omega$
(see  \cite[Theorem~3.3]{Fattorini-Russel}, \cite{Lebeau-Robbiano} and \cite{Fursikov-Imanuvilov-186}).
Degenerate parabolic equations exhibit a wider range of behaviours: observability  
may hold true or not depending on the type of degeneracy. For instance, 
the case of degenerate parabolic equations on the boundary of the domain in one space dimension 
is well understood (see \cite{Cannarsa-V-M-ADE}, \cite{Cannarsa-V-M-SIAM},  \cite{Ala-Can-Fra}, \cite{Martinez-Vancost-JEE-2006},  \cite{Can-Fra-Roc_1}, and \cite{Can-Fra-Roc_2}). Fewer results are available for multidimensional problems, see \cite{Cannarsa-V-M-CRAS}. 

As for parabolic equations with interior degeneracy, a fairly complete analysis is available for Grushin type operators 
\begin{equation} \label{Grushin}
\left\lbrace \begin{array}{ll}
\partial_t g - \Delta_x g - |x|^{2\gamma} \Delta_y g = 0 & \mbox{in}\; (0, \infty)\times\Omega\,,
\vspace{.1cm}\\
g(t,x,y)=0 & (t,x,y) \in (0, \infty)\times\partial \Omega\,,\vspace{.1cm}\\
g(0,x,y)=g^0(x,y)\,,                                                                                   & (x,y) \in \Omega\,,
\end{array}\right.
\end{equation}
where  
$\Omega:=\Omega_1 \times \Omega_2$, 
$\Omega_1$ is a bounded open subset of $\mathbb{R}^{N_1}$ such that $0 \in \Omega_1$,
$\Omega_2$ is a bounded open subset of $\mathbb{R}^{N_2}$, $N_1, N_2 \in \mathbb{N}^*:= \{1,2,3, ....\}$, and
$\gamma >0$. Indeed, it has been proved (\cite{Grushin, MR3162108}) that the observability inequality:
\begin{itemize}
\item holds in any positive time $T>0$ with an arbitrary open set $\omega\subset \Omega$  when $\gamma \in (0,1)$,
\item holds only in large time $T>T_{\min}>0$ when $\gamma=1$ and 
$\omega:= \omega_1 \times \Omega_2$ is a strip parallel to the $y$-axis  not containing the line segment $x=0$, and 
\item does not hold when $\gamma>1$.
\end{itemize}
Moreover, the  value of $T_{\min}$ has been explicitly computed for suitable observation regions $\omega$, see~\cite{minimal}.
The above observability properties may be changed by adding a zero order term with singular coefficient, see \cite{MR3205097} and  \cite{morgan}.
Similar results have been obtained for Kolmogorov type equations, see \cite{MR3163490, MR2566710, Kolmogorov}.

\subsubsection{Lipschitz stability}

Our formulation of the inverse problem corresponds to a single measurement
(see also Bukhgeim and Klibanov~\cite{BK} who first proposed a methodology 
based on Carleman estimates).  Following \cite{BK},  
many works have been published on this subject.  For uniformly parabolic equations we can refer the reader, for example, to 
Imanuvilov and Yamamoto~\cite{imayam98}, Isakov~\cite{Is},
Klibanov~\cite{Kl}, Yamamoto~\cite{Ya}, and the references therein (the present list of references is by no means complete).
Inverse problems for boundary-degenerate parabolic equations were studied by Cannarsa, Tort and Yamamoto~\cite{CTY1,CTY2}.
For Grushin type equations, the inverse source problem  was addressed in \cite{MR3162108}, and an inverse coefficient problem in \cite{inversecoefficient}.

\subsection{Structure of this article}
This paper has much to do with estimates. So, keeping track of all constants is definitely an issue. That is why we shall use capital letters, possibly with a subscript, only for those constants $C$ that are used in different parts of the article. Technical constants $c$ that are used in a single proof will be labeled by lower case letters, possibly with a subscript.

Sections \ref{se:wp} and \ref{sec:Fourier} are devoted to preliminary results concerning the well posedness of \eqref{H},
the Fourier decomposition of its solutions, and the dissipation speed of the Fourier modes.

In Section \ref{Sec:Heat_1D}, we state a Carleman estimate for a 1D-heat equation
with parameters $(n,p)$, solved by the Fourier modes of the solution of \eqref{H}.

In Section \ref{Sec:slice}, we prove Lipschitz stability with observation on a slice parallel to the $(y,z)$-plane (\ref{Thm:Lips_Stab_simple}).

In Section \ref{Sec:tube}, we prove Lipschitz stability with observation on a tube parallel to the $z$-axis (Theorem \ref{Thm:Lips_Stab}).

In Section \ref{Sec:Obs}, we prove that observability holds only in large time (Theorem~\ref{Thm:Obs}).

In Section \ref{sec:z_dans_R}, we state and justify analogous results for the Heisenberg heat equation  on $(-1,1)\times\mathbb{T}\times\mathbb{R}$.
Such a formulation allows to use the above theory to treat the Heisenberg equation written in the alternative form 
$$\Big(\partial_t  - \big(\partial_{x_1}-\frac{x_2}{2} \partial_{x_3} \big)^2 - \big(\partial_{x_2}+\frac{x_1}{2} \partial_{x_3} \big)^2 \Big) G(t,x_1,x_2,x_3) = 0.$$

Finally, in Section \ref{Sec:ccl}, we discuss conclusions an open problems.

%
%

\section{Well-posedness and unique continuation}\label{se:wp}
Without further specification, all functions are understood to be real-valued. 

\subsection{Well-posedness}

In this section, we recall well-posedness and regularity results for problem \eqref{H}. It is convenient to denote by $L^2([-1,1]\times \T\times\T)$, or briefly $L^2(\Omega)$, the space of all (equivalence classes of) Lebesgue-measurable functions $u:[-1,1]\times \R\times\R\to \R$ such that, for all $h,k\in\Z$,
\begin{equation}
\label{eq:per}
u(x,y+2h\pi,z+2k\pi)=u(x,y,z)\qquad (x,y,z)\in [-1,1]\times \R\times\R\;\;\mbox{a.e.}
\end{equation}
and
\begin{equation*}
\|u\|^2:=\int_{-1}^1dx\int_{-\pi}^\pi dy\int_{-\pi}^\pi |u(x,y,z)|^2dz<\infty.
\end{equation*}
$L^2(\Omega)$ is a Hilbert space over $\R$ with scalar product
\begin{equation*}
\langle u,v\rangle=\int_{-1}^1dx\int_{-\pi}^\pi dy\int_{-\pi}^\pi u(x,y,z)v(x,y,z)\,dz\qquad\forall \,u,v\in L^2(\Omega)\,.
\end{equation*}
Such a space will be also denoted by $\sph $. Now, consider the dense subspace $\mathcal C^\infty_{(0)}(\Omega)$ of $\sph $ which consists of all functions
$u\in\mathcal C^\infty\big([-1,1]\times \R\times\R\big)$ satisfying \eqref{eq:per} such that, for some $r\in [0,1)$,
\begin{equation*}
u(x,y,z)=0 \qquad\forall\,(x,y,z)\in \big([-1,1]\setminus [-r,r]\big)\times\R \times\R\,.
\end{equation*}
The bilinear form $ (\cdot,\cdot):\mathcal C^\infty_{(0)}(\Omega)\times  \mathcal C^\infty_{(0)}(\Omega)\to\R$ defined by
\begin{equation*}
(u,v)=\int_{-1}^1dx\int_{-\pi}^\pi dy\int_{-\pi}^\pi \big\{\partial_xu\partial_xv+(\partial_yu+x\partial_zu)(\partial_yv+x\partial_zv)\big\}\,dz
\end{equation*}
is positive definite because, for all $u\in \mathcal C^\infty_{(0)}(\Omega)$ we have
\begin{equation}
\label{eq:poincare}
\|u\|^2\leqslant 4 \int_{-1}^1dx\int_{-\pi}^\pi dy\int_{-\pi}^\pi|\partial_xu|^2dz\leqslant 4(u,u)\,.
\end{equation}
Denoting by $|\cdot|$ the norm associated with the scalar product $(\cdot,\cdot)$, we introduce the space $H^1_{(0)}(\Omega)$, or $\spv$,  as the closure of
$\mathcal C^\infty_{(0)}(\Omega)$ with respect to $|\cdot|$.  Observe that two bounded linear operators $X_1,X_2:\spv\to \sph $ are  defined by
\begin{equation*}
X_1u =\lim_{k\to\infty}\partial_xu_k\quad\mbox{and}\quad
X_2u =\lim_{k\to\infty}(\partial_yu_k+x\partial_zu_k)\,,
\end{equation*}
where $\{u_k\}_k$ is any sequence in $\mathcal C^\infty_{(0)}(\Omega)$ such that $|u_k-u|\to 0$ as $k\to\infty$. Moreover,   
\begin{equation*}
\int_{-1}^1dx\int_{-\pi}^\pi dy\int_{-\pi}^\pi  (X_1u )\,v\,dz
= -\int_{-1}^1dx\int_{-\pi}^\pi dy\int_{-\pi}^\pi u\,\partial_xv\,dz\,,
\end{equation*}
and
\begin{equation*}
\int_{-1}^1dx\int_{-\pi}^\pi dy\int_{-\pi}^\pi (X_2u)\,v\,dz
= -\int_{-1}^1dx\int_{-\pi}^\pi dy\int_{-\pi}^\pi u\,(\partial_yv+x\partial_zv)\,dz\,,
\end{equation*}
for all $u\in \spv$ and $v\in \mathcal C^\infty_{(0)}(\Omega)$. Also,
 the inequality 
\begin{equation*}
\|u\|\leqslant 2\|X_1u\|\qquad\forall \,u\in \spv
\end{equation*}
readily follows from \eqref{eq:poincare}. So, $\spv$ is a Hilbert space with  the scalar product  
\begin{equation*}
(u,v)=\int_{-1}^1dx\int_{-\pi}^\pi dy\int_{-\pi}^\pi \big\{(X_1u)(X_1v)+(X_2u)(X_2v)\big\}\,dz\qquad\forall \,u,v\in \spv\,.
\end{equation*}
Following a well-known procedure (\cite{lions61}) we can introduce the {\em regularly accretive operator} $\hl :D(\hl )\subset \sph \to \sph $ defined by
\begin{equation}\begin{cases}
\hspace{.cm}
D(\hl ) = \big\{u\in \spv~:~\exists\, C>0 \text{ such that } |(u,v)| \le C \|v\|\,,\;\;\forall v\in \spv\big\}
\\
\hl u=f \qquad\forall u\in D(\hl )\,,
\hspace{.cm}
\end{cases}
\end{equation}
where $f$ is the unique element of $\sph $ associated (via the Riesz isomorphism) with the extension to $\sph $ of the bounded linear functional 
$v\mapsto (u,v)$. Observe that $D(\hl )$ is dense in $\sph $ because it contains $C^\infty_{(0)}(\Omega)$. Therefore,
$\hl $ is a positive self-adjoint operator on $\sph $ satisfying $D(\hl ^{1/2})=\spv$ (\cite[Theorem~2.2.3]{tanabe}),
 and
  $-\hl $ generates an analytic semigroup of contractions on $\sph $ (\cite[Theorem~3.6.1]{tanabe}) that will be denoted by $S(t)$. 
  
 For every $g^0 \in \sph $ and $\widetilde{h} \in L^2(0,T;\sph )$, problem  \eqref{H} can be recast as follows
 \begin{equation}\label{eq:abstract_H}
\begin{cases}
\hspace{.cm}
g'(t)+\hl g(t)=\widetilde{h}(t)&t\in(0,T)
\\
\hspace{.cm}
g(0)=g^0\,.
\end{cases}
\end{equation}
The function $g \in C^0([0,T];\sph ) \cap L^2(0,T;\spv)$ given by
\begin{equation*}
g(t)=S(t)g^0+\int_0^tS(t-s)\widetilde{h}(s)ds\qquad t\in[0,T]
\end{equation*}
is called the {\em mild solution} of \eqref{eq:abstract_H}.
It is well known that the mild solution of \eqref{eq:abstract_H} is also a {\em weak solution} in the following sense:
for every $v\in D(\hl )$
\begin{itemize}
\item  the function $\langle g(\cdot),v\rangle$ is absolutely continuous on $[0,T]$, and 
\item for a.e. $t\in [0,T]$
\begin{equation}\label{eq:defweaksol}
\frac{d}{dt} \langle g(t),v\rangle +\langle g(t),\hl  v\rangle = \langle \widetilde{h}(t),v\rangle\, .
\end{equation}
\end{itemize}
Note that, as showed in \cite{lions61}, condition \eqref{eq:defweaksol} is equivalent to the definition of solution by transposition, that is,
$$
\begin{array}{l}
\displaystyle\int_{-1}^1dx\int_{-\pi}^\pi dy\int_{-\pi}^\pi 
\Big\{ g(\tau,x,y,z)\varphi(\tau,x,y,z) - g^0(x,y,z) \varphi(0,x,y,z) \Big\}dz \\
\displaystyle  =\int_0^{\tau}dt\int_{-1}^1dx\int_{-\pi}^\pi dy\int_{-\pi}^\pi 
g\Big\{\partial_x^2\varphi - (x \partial_z+\partial_y)^2\varphi \Big\} \,dz
\end{array}
$$
for every  $\tau \in (0,T)$ and every function 
$\varphi \in \mathcal C^{2}\big([0,T] \times  [-1,1]\times \T\times\T\big)$.

The following proposition describes  well-known properties of  mild solutions that follow from the analiticity of $S(t)$.
\begin{Prop} \label{Prop:existence}
For every $g^0 \in \sph $, $T>0$, and 
$\widetilde{h} \in L^{2}(0,T;\sph )$,
the mild solution $g$ of the Cauchy problem \eqref{eq:abstract_H} satisfies
\begin{equation}\label{C0continforsol}
\|g(t)\|
\leqslant   \|f_0\| + \sqrt{T} \|\widetilde{h}\|_{L^2(0,T;\sph )} \quad\forall t \in [0,T]\,.
\end{equation}
Moreover, for every $\tau\in(0,T]$,
\begin{equation*}
g\in H^1(\tau,T;\sph )\cap \mathcal C([\tau,T];\spv) \cap L^2(\tau,T;D(\hl ))\,.
\end{equation*}
In particular,
$g(t) \in D(\hl )$ and $g'(t) \in \sph $ for a.e. $t\in [0,T]$. 
\end{Prop}

\subsection{Unique continuation}

Observe that, in particular, \eqref{eq:observability} yields a unique continuation property for \eqref{H}.
The following more general result, which is a consequence of Holmgren's uniqueness theorem, suggests that no obstruction to 
observability  should be expected for problem \eqref{H}. The proof is given in the appendix,  Section \ref{Appendix:UC}.

\begin{Prop} \label{Prop:Holmgren}
Let $T>0$  and let $\omega$ be as in \eqref{eq:omega}.
Any solution $$g \in C^0([0,T],L^2(\Omega))$$ of \eqref{H} with $\tilde{h}=0$,
which vanishes on $(0,T)\times\omega$ is identically zero.
\end{Prop}

\section{Fourier decomposition and dissipation} \label{sec:Fourier}

\subsection{Fourier decomposition}

We are now going to study the Fourier decomposition of the solution of \eqref{H}. For this purpose, for any $(n,p) \in \mathbb{Z}^2$ let us consider the operator $$\hl _{n,p}:D(\hl _{n,p})\subset L^2(-1,1;\C)\to L^2(-1,1;\C)$$ defined by
\begin{equation}\begin{cases}
\hspace{.cm}
D(\hl _{n,p}) = H^2\cap H^1_0(-1,1;\C)
\\
\hl _{n,p}u(x)=- u''(x) + (px+n)^2u(x) \qquad\forall u\in D(\hl _{n,p})\,.
\hspace{.cm}
\end{cases}
\end{equation}
It is well known that $\hl _{n,p}$ is a positive self-adjoint operator on $L^2(-1,1;\C)$ 
 and
  $-\hl _{n,p}$ generates an analytic semigroup of contractions. The notion of mild/weak solutions of the evolution equation associated
  with $\hl _{n,p}$, that we recalled in section~\ref{se:wp}, is used in our next proposition.

\begin{Prop} \label{Prop:Fourier}
Let $g^0 \in \sph =L^2(\Omega)$, $T>0$, and 
$\widetilde{h} \in L^{2}\big((0,T)\times \Omega\big)$. Then
the mild solution $g$ of the Cauchy problem \eqref{H} satisfies, in $L^{2}\big((0,T)\times \Omega\big)$,
\begin{equation}\label{eq:Fourier_g}
g(t,x,y,z)=  \sum\limits_{n,p \in \mathbb{Z}} g_{n,p}(t,x) e^{i(ny+pz)}
\end{equation}
where
\begin{equation}\label{eq:Fourier_coeff}
g_{n,p}(t,x):=\frac{1}{(2\pi)^2} \int_{\mathbb{T}^2} g(t,x,y,z) e^{-i(ny+pz)} dy dz
\end{equation}
belongs to $\mathcal C\big([0,T];L^2(-1,1;\C)\big)\cap L^2\big(0,T;H^1_0(-1,1;\C)\big)$.
Moreover, for every $(n,p) \in \mathbb{Z}^2$, $g_{n,p}$ is the mild solution of the Cauchy problem
\begin{equation} \label{H_np}
\left\lbrace\begin{array}{ll}
\Big(\partial_t  - \partial_x^2 + (px+n)^2 \Big) g_{n,p}(t,x) = \widetilde{h}_{n,p}(t,x)\,, & (t,x) \in (0,T)\times(-1,1) \,, \\
g_{n,p}(t,\pm 1)=0\,,                                                                       & t \in (0,T)\,, \\
g_{n,p}(0,x)=g_{n,p}^0(x)\,,                                                                & x \in (-1,1)\,,
\end{array}\right.
\end{equation}
where
\begin{equation*}
\widetilde{h}_{n,p}(t,x):=\frac{1}{(2\pi)^2} \int_{\mathbb{T}^2} \widetilde{h}(t,x,y,z) e^{-i(ny+pz)} dy dz
\end{equation*}
and
\begin{equation*}
g^0_{n,p}(x):=\frac{1}{(2\pi)^2} \int_{\mathbb{T}^2} g^0(x,y,z) e^{-i(ny+pz)} dy dz\,.
\end{equation*}
Furthermore, if 
\begin{equation*}
g^0_{n,p}\in H^2\cap H^1_0(-1,1;\C)\quad\mbox{and}\quad
\widetilde{h}_{n,p}\in H^1\big(0,T;L^2(-1,1;\C)\big),
\end{equation*}
then the function $v_{n,p}:=\partial_t g_{n,p}\in \mathcal C\big([0,T];L^2(-1,1;\C)\big)\cap L^2\big(0,T;H^1_0(-1,1;\C)\big)$ is the weak solution of
\begin{equation*} 
\left\lbrace\begin{array}{ll}
\Big(\partial_t  - \partial_x^2 + (px+n)^2 \Big) v_{n,p}(t,x) = \partial_t\widetilde{h}_{n,p}(t,x)\,, & (t,x) \in (0,T)\times(-1,1) \,, \\
v_{n,p}(t,\pm 1)=0\,,                                                                       & t \in (0,T)\,, \\
v_{n,p}(0,x)=\hl _{n,p}g_{n,p}^0(x)+\widetilde{h}_{n,p}(0,x)\,,                                                                & x \in (-1,1)\,.
\end{array}\right.
\end{equation*}
\end{Prop}

\noindent \textbf{Proof:} The relations \eqref{eq:Fourier_g} and \eqref{eq:Fourier_coeff} between $g$ and the family of the Fourier coefficients of $g(t,x,\cdot,\cdot)$ is justified by the fact that $g\in C\big([0,T];L^2(\Omega)\big)$. Thus, we just have to show that, for every $(n,p) \in \mathbb{Z}^2$, $g_{n,p}$ is the weak solution of  problem \eqref{H_np}. Since
\begin{equation*}
g_{n,p}\in \mathcal C\big([0,T];L^2(-1,1;\C)\big)\cap L^2\big(0,T;H^1_0(-1,1;\C)\big)
\end{equation*}
in view of \eqref{eq:Fourier_coeff}, we just need to show that, for every $\varphi\in H^2\cap H^1_0(-1,1;\C)$,
\begin{itemize}
\item[(i)]  the function $t\mapsto\int_{-1}^1 g_{n,p}(t,x)\varphi(x)\,dx$ is absolutely continuous on $[0,T]$,  
\item[(ii)] for a.e. $t\in [0,T]$
\begin{equation*}
\frac{d}{dt} \int_{-1}^1 g_{n,p}(t,x)\varphi(x)\,dx + \int_{-1}^1 g_{n,p}(t,x)\hl _{n,p}\varphi(x)\,dx = \int_{-1}^1 \widetilde{h}(t,x)\varphi(x)\,dx \, .
\end{equation*}
\end{itemize}
Indeed, since
\begin{equation*}
\int_{-1}^1 g_{n,p}(t,x)\varphi(x)\,dx=
\frac{1}{(2\pi)^2} \int_{-1}^1dx\int_{\mathbb{T}^2} g(t,x,y,z) e^{-i(ny+pz)}\varphi(x)\,dx dy dz\,,
\end{equation*}
property (i) follows from the fact that $g$ is the weak solution of \eqref{eq:abstract_H} and the real and imaginary parts, $u$ and $v$, of the complex-valued function
\begin{equation*}
w(x,y,z):=e^{-i(ny+pz)}\varphi(x)\qquad (x,y,x)\in\Omega
\end{equation*}
belong to $D(\hl )$. As for property (ii), observe that by the same argument
\begin{eqnarray*}
\lefteqn{\frac{d}{dt} \int_{-1}^1 g_{n,p}(t,x)\varphi(x)\,dx + \int_{-1}^1 g_{n,p}(t,x)\hl _{n,p}\varphi(x)\,dx}
\\
&=&\frac{1}{(2\pi)^2} \frac{d}{dt} \int_{-1}^1dx\int_{\mathbb{T}^2} g(t,x,y,z) e^{-i(ny+pz)}\varphi(x)\,dy dz
\\
&&+\frac{1}{(2\pi)^2} \int_{-1}^1dx\int_{\mathbb{T}^2} g(t,x,y,z) e^{-i(ny+pz)}\big(- \varphi''(x) + (px+n)^2\varphi(x)\big)\,dy dz
\\
&=&\frac{1}{(2\pi)^2} \frac{d}{dt} \int_{-1}^1dx\int_{\mathbb{T}^2} g(t,x,y,z) \big(u(x,y,z)+iv(x,y,z)\big)\,dy dz
\\
&&+\frac{1}{(2\pi)^2} \int_{-1}^1dx\int_{\mathbb{T}^2} g(t,x,y,z) \big(\hl u(x,y,z)+i\hl v(x,y,z)\big)\,dy dz
\\
&=&\frac{1}{(2\pi)^2}\int_{-1}^1dx\int_{\mathbb{T}^2}  \widetilde{h}(t,x,y,z) e^{-i(ny+pz)}\varphi(x)\,dy dz
\\
&=& \int_{-1}^1  \widetilde{h}_{n,p}(t,x)\varphi(x)\,dx
\end{eqnarray*}
for a.e. $t\in [0,T]$. This completes the proof.  \hfill $\Box$

\subsection{Dissipation speed on $(-1,1)$}

For any $(n,p) \in \mathbb{Z} \times \mathbb{Z}$, we define
\begin{equation} \label{def:lambda}
\lambda_{n,p}=\inf_{\varphi \in H^1_0(-1,1)}\Big\{ 
\int\limits_{-1}^{1} \Big[ \varphi'(x)^2 + (px+n)^2 \varphi(x)^2 \Big] dx~:~\int_{-1}^1 \varphi(x)^2 dx =1
\Big\}.
\end{equation}

\begin{Prop} \label{Prop:dissip}
The following inequalities hold:
\begin{equation} \label{dissip_speed_p}
\lambda_{n,p} \geqslant \frac{1}{4} (|p|+1)\,, \quad \forall (n,p) \in \mathbb{Z} \times \mathbb{Z} 
\end{equation}
\begin{equation} \label{dissip_speed_n}
\lambda_{n,p} \geqslant \frac{n^2}{4} \,, \quad \forall (n,p) \in \mathbb{Z} \times \mathbb{Z} \text{ such that } |n| \geqslant 2 |p|.
\end{equation}
\end{Prop}

\begin{rk} \label{rk:validite_LR}\rm
Observe that, when $|n| \geqslant 2|p|$,  the dependence of $\lambda_{n,p}$ is quadratic with respect to $n$.
This is the key point to apply the Lebeau-Robbiano strategy with respect to the variable $y$ 
($n$ has to be negligeble with respect to $\lambda_{n,p}$ when $n \rightarrow \infty$ and $p$ is fixed).
This is no longer true when $x$ is free to range in the whole space $\mathbb{R}$ because of translation invariance.
\end{rk}

\noindent \textbf{Proof of Proposition \ref{Prop:dissip}:} If $p=0$ then
$\lambda_{n,0} \geqslant \big(\frac{\pi}{2})^2 + n^2$. So,  (\ref{dissip_speed_p}) and (\ref{dissip_speed_n}) hold true.
Let now $(n,p) \in \mathbb{Z} \times [\mathbb{Z} \setminus \{0\}]$ and observe that, without loss of generality, one may assume $p>0$.
By the change of variable 
$$\varphi(x)=\sqrt[4]{p} \;\psi\left( \tilde{x}=\sqrt{p}\Big(x+\frac{n}{p}\Big) \right),$$
from (\ref{def:lambda})   we deduce that
\begin{eqnarray*}
\lambda_{n,p} 
& 
\geqslant& \inf_{\varphi \in H^1(\mathbb{R}) \cap L^2(\mathbb{R},x^2 dx)}\Big\{ 
\int\limits_{\mathbb{R}} \Big[ |\varphi'|^2 + (px+n)^2 |\varphi|^2 \Big] dx~:~\int\limits_{\mathbb{R}} |\varphi|^2 =1 \Big\}
\\ 
& =& p \,  \inf_{\psi \in H^1(\mathbb{R}) \cap L^2(\mathbb{R},\tilde{x}^2 d\tilde{x})}
\Big\{ 
\int\limits_{\mathbb{R}} \big[ \psi'(\tilde{x})^2 + \tilde{x}^2 \psi(\tilde{x})^2 \big] d\tilde{x}~:~
\int\limits_{\mathbb{R}} |\psi|^2  = 1 \Big\},
\end{eqnarray*}
where we have denoted by $L^2(\mathbb{R},x^2 dx)$ the space of all Lebesgue measurable functions $\varphi:\R\to\R$ such that
\begin{equation*}
\int\limits_{\mathbb{R}}|\varphi(x)|^2 x^2 dx<\infty.
\end{equation*}
Since the last infimum above equals $1$, (\ref{dissip_speed_p}) is proved.
Now, suppose $|n| \geqslant 2p$. Then for every $x \in [-1,1]$
$$|px+n| \geqslant |n| - p \geqslant \frac{n}{2}\,.$$
Thus
$$\lambda_{n,p}  \geqslant \inf_{\varphi \in H^1_0(-1,1)}\Big\{ 
\frac{n^2}{4} \int\limits_{-1}^{1}  \varphi(x)^2 dx~:~\int_{-1}^1 \varphi(x)^2 dx =1
\Big\}
 = \frac{n^2}{4}\,,$$ 
which proves (\ref{dissip_speed_n}). \hfill $\Box$

\section{1D heat equations with parameters}
\label{Sec:Heat_1D}

In this section, we will prove several estimates for 1D heat equations with parameters which
 will be used in the proof of the main results of the paper.

\subsection{Carleman estimates}

Let us set $\R_+=(0,\infty)$. For a given $T>0$ and any $(n,p) \in \mathbb{Z} \times \mathbb{R}_+$, we define the operator
$$\mathcal{P}_{n,p} g = \partial_t g - \partial_x^2 g + (px+n)^2 g$$
acting on functions $g: [0,T]\times [-1,1]\to\C$.
\begin{Prop} \label{Carleman_global}
Let $a, b\in\R$ be such that $-1 \leqslant a < b  \leqslant 1$.
Then there exist a weight function $\beta \in \mathcal C^3([-1,1];\mathbb{R}_+)$ and
positive constants $C_1, C_2$ such that
for any $(n,p) \in \mathbb{Z} \times \mathbb{R}_+$, any $T>0$, and any 
$$g \in \mathcal C([0,T];L^2(-1,1)) \cap L^2(0,T;H^1_0(-1,1))$$
the following inequality holds
\begin{eqnarray} \label{Carl_est}
\lefteqn{C_1 \int_0^T\!\!\! \int_{-1}^1 \Big(
\frac{M}{t(T-t)} \,\left|\partial_xg\right|^2  +
\frac{M^3}{[t(T-t)]^3} \left|g\right|^2 
\Big)  e^{-\frac{M \beta(x)}{t(T-t)}} dx dt}
\\
& & \leqslant
\int_0^T\!\!\! \int_{-1}^1 | \mathcal{P}_{n,p} g |^2 e^{-\frac{ M \beta(x)}{t(T-t)}} dx dt +
\int_0^T\!\!\! \int_a^b \frac{M^3}{[t(T-t)]^3} |g|^2 e^{-\frac{M \beta(x)}{t(T-t)}} dx dt
\nonumber
\end{eqnarray}
where 
\begin{equation} \label{M_Carleman}
M:=C_2 \max\{T+T^2;(|n|+p)T^2\}\,.
\end{equation}
\end{Prop}
In the appendix, we give a complete proof of the above Carleman estimate.
\begin{rk}\rm
The proof of the main results of this article
only uses the above result  for $p \in \mathbb{Z}$.
However, we prefer to derive most of our preliminary results for $p \in \mathbb{R}$ instead of $p \in \mathbb{Z}$ in order to
justify the generalization discussed in Section \ref{sec:z_dans_R},
where the domain is $(-1,1)\times\mathbb{T}\times\mathbb{R}$.
\end{rk}

\subsection{1D observability inequality with source term}

The goal of this section is the proof of the following result.

\begin{Prop} \label{Prop:1D-IO-source}
Let $a, b\in\R$ be such that $-1 \leqslant a < b  \leqslant 1$.
Then there exist constants $C_3, C_4 >0$ such that,
for every $T>0$, $p \in \mathbb{R}$, $n \in \mathbb{Z}$, $g_{n,p}^0 \in L^2(-1,1)$, and $\widetilde{h}_{n,p} \in L^2\big([0,T]\times[-1,1]\big)$
the solution of $(\ref{H_np})$ satisfies
\begin{eqnarray} \label{1D-IO-source}
\nonumber
\int\limits_{-1}^1 |g_{n,p}(T,x)|^2 dx \leqslant 
 e^{C_3 \left( 1 + \frac{1}{T} + |p|-C_4\min\{|p|\,,\,p^2\} T \right)} \int\limits_0^T \!\int\limits_a^b |g_{n,p}(t,x)|^2 dx dt
\\
+ \epsilon_{n,p}(T) \int\limits_0^T \!\int\limits_{-1}^1 |\widetilde{h}_{n,p}(t,x)|^2 dx dt\,,
\end{eqnarray}
for some constant $\epsilon_{n,p}(T)$ satisfying
\begin{equation} \label{def:epsilon_np_n<2p}
|\epsilon_{n,p}(T)| \leqslant  \frac{C_3}{|p|+1} + 
e^{C_3 \left( 1 + \frac{1}{T} + |p|-C_4\min\{|p|\,,\,p^2\} T \right)} =: \epsilon_p'(T)
\end{equation}
for all $(n,p) \in \mathbb{Z}\times\R$ and
\begin{equation} \label{def:epsilon_np_n>2p}
|\epsilon_{n,p}(T)| \leqslant  \frac {C_3}{n^2} + e^{C_3 \left( 1 + \frac{1}{T} - C_4 n^2 T \right)} =:\epsilon_n''(T) \;\text{ if }\; |n| > 2|p|\,.
\end{equation}
\end{Prop}

We will use the following preliminary result.

\begin{Lem}\label{Lem:Duhamel_np}
For every $0 \leqslant T_1 < T_2 < \infty$, $(n,p) \in \mathbb{Z}\times\R$, $g_{n,p}^0 \in L^2(-1,1)$, and $\widetilde{h}_{n,p} \in L^2\big([0,T]\times[-1,1]\big)$
the solution of $(\ref{H_np})$ satisfies
$$\|g_{n,p}(T_2)\|^2 \leqslant 2 \|g_{n,p}(T_1)\|^2 e^{-2\lambda_{n,p}(T_2-T_1)} 
+ \frac{1}{\lambda_{n,p}} \int_{T_1}^{T_2} \|\widetilde{h}_{n,p}(t)\|^2 dt$$
where $\|.\|=\|.\|_{L^2(-1,1)}$.
\end{Lem}

\noindent \textbf{Proof of Lemma \ref{Lem:Duhamel_np}:} By Duhamel's formula and the Cauchy-Schwarz inequality we obtain 
$$\begin{array}{ll}
\|g_{n,p}(T_2)\|
& \leqslant e^{-\lambda_{n,p}(T_2-T_1)} \|g_{n,p}(T_1)\| + \int_{T_1}^{T_2} e^{-\lambda_{n,p}(T_2-t)} \|\widetilde{h}_{n,p}(t)\| dt
\\ & \leqslant
e^{-\lambda_{n,p}(T_2-T_1)} \|g_{n,p}(T_1)\|  + \frac{1}{\sqrt{2\lambda_{n,p}}} \left( \int_{T_1}^{T_2} \|\widetilde{h}_{n,p}(t)\|^2 dt \right)^{1/2}.
\end{array}$$
The inequality $(a+b)^2 \leqslant 2 a^2 + 2 b^2$ gives the conclusion. $\hfill \Box$
\\

\noindent \textbf{Proof of Proposition \ref{Prop:1D-IO-source}:} 
In this proof, we write $g$, $g^0$, $\widetilde{h}$ and $\|.\|$ 
instead of $g_{n,p}$, $g_{n,p}^0$, $\widetilde{h}_{n,p}$ and $\|.\|_{L^1(-1,1)}$ in order to simplify the notation.
\\

\noindent \emph{Step 1: use of dissipation.}
Appying Lemma \ref{Lem:Duhamel_np} with $(T_1,T_2)=(t,T)$
and integrating the resulting inequality over $t \in (T/3,2T/3)$ yields 
\begin{equation} \label{step1}
\int\limits_{-1}^1 |g(T)|^2 dx \leqslant 
\frac{6}{T}  e^{-2 \lambda_{n,p} \frac{T}{3}} \int\limits_{T/3}^{2T/3}\! \int\limits_{-1}^1 |g|^2 dx dt  + 
\frac{1}{\lambda_{n,p}} \int\limits_0^T\!\int\limits_{-1}^1 |\widetilde{h}|^2 dx dt\,.
\end{equation}

\noindent \emph{Step 2: existence of a constant $c_0>0$, independent of $(T,n,p,g^0,\widetilde{h})$, such that}
\begin{equation} \label{step2}
\int\limits_{T/3}^{2T/3} \!\int\limits_{-1}^1 |g|^2 dx dt 
\leqslant  c_0 T e^{\frac{9M\beta^*}{2T^2}} \Big( 
\int\limits_0^T\! \int\limits_a^b |g|^2 dx dt + \int\limits_0^T\! \int\limits_{-1}^1  |\widetilde{h}|^2 dx dt  \Big)
\end{equation}
where $\beta$, $C_2$, and $M=M(n,p,T)$ are as in Proposition \ref{Carleman_global} and 
$$\beta^*:=\max\{\beta(x);x \in [-1,1]\}\,.$$
By Proposition \ref{Carleman_global} we have
\begin{eqnarray*}
\lefteqn{ C_1 \Big( \frac{4M}{T^2} \Big)^3 e^{-\frac{9M\beta^*}{2T^2}} \int_{T/3}^{2T/3}\!\!\! \int_{-1}^1 |g|^2 dxdt 
 \leqslant
C_1 \int_{T/3}^{2T/3}\!\!\! \int_{-1}^1 \frac{M^3}{[t(T-t)]^3} |g|^2  e^{-\frac{M \beta(x)}{t(T-t)}} dxdt}
\\ &\leqslant &
C_1 \int^T_0 \!\!\!  \int_{-1}^1 \frac{M^3}{[t(T-t)]^3} |g|^2 e^{-\frac{M\beta(x)}{t(T-t)}} dxdt 
\\ 
&\leqslant &\int_0^T\!\!\!  \int_{-1}^1 |\widetilde{h}|^2 e^{-\frac{M \beta(x)}{t(T-t)}}  dxdt 
+ \int_0^T \!\!\! \int_a^b \frac{M^3}{[t(T-t)]^3} |g|^2  e^{-\frac{M \beta(x)}{t(T-t)}} dxdt 
\\ 
&\leqslant & \int_0^T\!\!\!  \int_{-1}^1 |\widetilde{h}|^2 dxdt + c_1\int_0^T\!\!\!  \int_a^b |g|^2 dx dt
\end{eqnarray*}
where $c_1=\sup\{ x^3 e^{-\beta_* x} ; x \geqslant 0\}$ and $\beta_*:=\min\{ \beta(x);x \in (a,b)\}$. Thus
$$\int_{T/3}^{2T/3} \!\!\int_{-1}^1 |g(t)|^2 dxdt 
\leqslant \frac{\max\{1,c_1\}}{4^3 C_1} \frac{T^6}{M^3} e^{\frac{9M\beta^*}{2T^2}} \Big( 
\int_0^T\!\!\!  \int_{-1}^1 |\widetilde{h}|^2  + \int_0^T\!\!\!  \int_a^b |g|^2  \Big).$$
We remark that $M \geqslant C_2 T$ and $M \geqslant C_2 T^2$ thus 
$T^6/M^3 \leqslant T/C_2^3$. Then, the previous inequality
gives (\ref{step2}) with $c_0=\max\{1,c_1\}/(4^3 C_1 C_2^3)$.
\\

\noindent \emph{Step 3: combination of \eqref{step1} and \eqref{step2}.} 
\begin{equation} \label{Obs_interm}
\begin{array}{ll}
\displaystyle \int_{-1}^1 |g(T)|^2 dx  \leqslant 
& 
\displaystyle 6 \,c_0\,  e^{\frac{9M\beta^*}{2T^2}-2 \lambda_{n,p} \frac{T}{3}} \int_0^T \int_a^b |g|^2 dxdt
\\  &
\displaystyle + \Big( \frac{1}{\lambda_{n,p}} + 6 \,c_0\,  e^{ \frac{9M\beta^*}{2T^2}-2 \lambda_{n,p} \frac{T}{3} }  \Big) 
\int_0^T \int_{-1}^1 |\widetilde{h}|^2 dx dt\,.
\end{array}
\end{equation}
From now on, we introduce the constants
$$C_3:=\ln(6c_0)+\frac{27 C_2 \beta^*}{2} + 3 \alpha^2 + 4  \quad \text{ and } \quad C_4:=\frac{1}{12 C_3}$$
where $C_2$ is as in (\ref{M_Carleman}) and $\alpha:=27 \beta^* C_2/4$.
\\

\noindent \emph{Step 4: proof of}
\begin{equation} \label{Majo_avec_C3C4_1}
\ln(6c_0)+\frac{9M\beta^*}{2T^2}-2 \lambda_{n,p} \frac{T}{3} \leqslant 
C_3 \left( 1 + \frac{1}{T} + |p|-C_4\min\{|p|\,,\,p^2\} T \right)\,, 
\end{equation}
for every $(n,p) \in \mathbb{Z}\times\R$.
\\

\noindent \emph{\underline{Case 1: $|n|< 2|p|$}}. By (\ref{dissip_speed_p}) and (\ref{M_Carleman}) we have that
$$M(T,n,p) \leqslant C_2 \left( T^2 + T + 3 |p| T^2 \right)
\quad \text{ and } \quad \lambda_{n,p} \geqslant  \frac{1}{4}(|p|+1)\,.$$
Thus 
\begin{eqnarray*}
\ln(6\,c_0) + \frac{9M\beta^*}{2T^2}- 2 \lambda_{n,p}\frac{T}{3} 
 \leqslant \ln(6\,c_0) + \frac{9 C_2 \beta^*}{2} \Big( 1 + \frac{1}{T} + 3 |p| \Big) - \frac{1}{6} |p| T
 \\
  \leqslant  C_3 \Big( 1+\frac{1}{T} + |p| - C_4 |p| T\Big) 
\end{eqnarray*}
which gives (\ref{Majo_avec_C3C4_1}).
\\

\noindent \emph{\underline{Case 2: $|n|\geqslant 2|p|$}}. In view of (\ref{dissip_speed_n}) and (\ref{M_Carleman}),
$$M(T,n,p) \leqslant C_2 \Big( T^2 + T + \frac{3}{2} |n| T^2 \Big)
\quad \text{ and } \quad \lambda_{n,p} \geqslant \frac{n^2}{4}\,.$$
Therefore,
\begin{eqnarray} \label{etoile}
\nonumber
\ln(6\,c_0) + \frac{9M\beta^*}{2T^2}- 2 \lambda_{n,p}\frac{T}{3} 
 \leqslant 
  \ln(6\,c_0) + \frac{9 C_2 \beta^*}{2} \Big( 1 + \frac{1}{T} + \frac{3}{2} |n| \Big) - \frac{1}{6}  n^2 T
 \\
  \leqslant \ln(6\,c_0) + \frac{9 C_2 \beta^*}{2} \Big( 1 + \frac{1}{T} \Big) + \frac{3 \alpha^2}{T} - \frac{1}{12}  n^2 T
\end{eqnarray}
because the maximal value of the function $f:s \in (0,\infty) \mapsto \alpha s - \frac{s^2 T}{12}$ is exactly $\frac{3 \alpha^2}{T}$.
Finally, using the assumption $|n|\geqslant 2|p|$, we obtain
\begin{eqnarray*}
\ln(6\,c_0) + \frac{9M\beta^*}{2T^2}- 2 \lambda_{n,p}\frac{T}{3} 
\leqslant
\ln(6\,c_0) + \frac{9 C_2 \beta^*}{2} \Big( 1 + \frac{1}{T} \Big) + \frac{3 \alpha^2}{T} - \frac{1}{3}  p^2 T 
\\
\leqslant C_3 \left( 1 + \frac{1}{T} - C_4 p^2 T \right)\,,
\end{eqnarray*}
which gives (\ref{Majo_avec_C3C4_1}).
\\

\noindent \emph{Step 5: proof of}
\begin{equation} \label{Majo_avec_C3C4_1}
\ln(6c_0)+\frac{9M\beta^*}{2T^2}-2 \lambda_{n,p} \frac{T}{3} \leqslant 
C_3 \left( 1 + \frac{1}{T} - C_4 n^2 T \right)\,, 
\end{equation}
for every $(n,p) \in \mathbb{Z}\times\R$ with $|n| \geqslant 2|p|$.
It results from (\ref{etoile}) and the choice of $C_3$ and $C_4$.
\\

\noindent \emph{Step 6: conclusion.}
From (\ref{Obs_interm}), Step 4 and (\ref{dissip_speed_p}), we deduce that (\ref{1D-IO-source}) and (\ref{def:epsilon_np_n<2p}) hold.
From (\ref{Obs_interm}), Step 4, Step 5 and (\ref{dissip_speed_n}), we obtain (\ref{1D-IO-source}) and (\ref{def:epsilon_np_n>2p}). $\hfill \Box$

\section{3D-Lipschitz stability estimate when $\omega$ is a slice}
\label{Sec:slice}

The goal of this section is the proof of Theorem \ref{Thm:Lips_Stab_simple}.
We focus on the uniform Lipschitz stability estimate for systems (\ref{H_np}) in the sense of the following definition.
We assume the source term $\widetilde{h}_{n,p}$ in (\ref{H_np}) takes the form
\begin{equation} \label{form_hnp}
\begin{array}{c}
\widetilde{h}_{n,p}(t,x)=R(t,x) h_{n,p}(x) 
\vspace{.2cm}
\\
\text{ where } h_{n,p} \in L^2(-1,1) \text{ and } R \in C^0([0,T] \times [-1,1])\,.
\end{array}
\end{equation}

\begin{Def}[Uniform Lipschitz stability]
Let $a,b \in \mathbb{R}$ with $-1 \leqslant a < b \leqslant 1$, $T>0$ and $0< T_0 < T_1 \leqslant T$.
We say the system \eqref{H_np}-\eqref{form_hnp} satisfies a \emph{uniform Lipschitz stability estimate on $(T_0,T_1)\times(a,b)$} if,
there exists a constant $C>0$ such that,
for every $p \in \mathbb{R}$, $n \in \mathbb{Z}$, $g_{n,p}^0 \in L^2(-1,1)$, $h_{n,p} \in L^2(-1,1)$, 
the solution of \eqref{H_np}-\eqref{form_hnp} satisfies
\begin{equation} \label{Lipschitz-ineq-n}
\int\limits_{-1}^1 |h_{n,p}|^2 dx 
\leqslant
C \Big(
\int\limits_{T_0}^{T_1} \int\limits_a^b |\partial_t g_{n,p}|^2 dx dt +
\int\limits_{-1}^1 | \hl _{n,p} g_{n,p}(T_1,x)|^2 dx 
\Big)
\end{equation}
where $\hl _{n,p}:=-\partial_x^2 +(px+n)^2$.
\end{Def}
Theorem \ref{Thm:Lips_Stab_simple} is a consequence of the next result and Bessel-Parseval identity.

\begin{Prop} \label{Prop:UO}
Let $a,b \in \mathbb{R}$ be such that $-1 < a < b < 1$ and $R$ be such that \eqref{Hyp:R} holds.
The exists $T^*>0$ such that, for every $T_0 \in (0,T_1-T^*)$
system \eqref{H_np}-\eqref{form_hnp} satisfies a uniform Lipschitz stability estimate on $(T_0,T_1) \times (a,b)$.
\end{Prop}
\begin{rk}\rm
 Inequality \eqref{Lipschitz-ineq-n}, with a constant $C$ that may depend on $n$ and $p$ is already known (see~\cite{imayam98}).
 Therefore, in order to prove Proposition~\ref{Prop:UO} it suffices to focus on high frequencies $(n,p)$.
\end{rk}

\noindent \textbf{Proof of Proposition \ref{Prop:UO}:} Let $C_3$ and $C_4$ be the constants given by Proposition~\ref{Prop:1D-IO-source}. We assume $(T_1-T_0)>T^*:=1/C_4$. 
\\

\noindent \emph{Step 1: application of Proposition~\ref{Prop:1D-IO-source}.}
From \eqref{Hyp:R} it follows that
$$R_0 |h_{n,p}(x)|  \leqslant |R(T_1,x) h_{n,p}(x)|  \leqslant |\partial_t g_{n,p}(T_1,x)| + |\hl _{n,p} g_{n,p}(T_1,x)|$$
and
\begin{equation} \label{interm1}
\int\limits_{-1}^1 |h_{n,p}|^2 dx 
\leqslant \frac{2}{R_0^2} \Big(
\int\limits_{-1}^1 |\partial_t g_{n,p}(T_1,x)|^2 dx +
\int\limits_{-1}^1 |\hl _{n,p} g_{n,p}(T_1,x)|^2 dx 
\Big).
\end{equation}
Notice that
$$
|p|-C_4\min\{\,|p|\,,\,p^2\} (T_1-T_0) \leqslant
\left\lbrace \begin{array}{ll}
1-C_4 p^2 (T_1 - T_0 - T_*) \leqslant 1 \, \text{ if } |p| \leqslant 1 
\vspace{.2cm}\\
-C_4 |p| (T_1-T_0-T_*) \leqslant 0\, \text{ if } |p| \geqslant 1\,.
\end{array}\right.$$
Thus, by Proposition~\ref{Prop:Fourier} and Proposition~\ref{Prop:1D-IO-source}, applied to $\partial_t g_{n,p}$,
we get
\begin{equation} \label{interm}
\begin{array}{ll}
\int\limits_{-1}^1 |\partial_t g_{n,p}(T_1,x)|^2 dx \leqslant 
& e^{C_3 \left( 2 + \frac{1}{T_1-T_0} \right)} \int\limits_{T_0}^{T_1} \int\limits_a^b |\partial_t g_{n,p}|^2 dx dt \\
&  + \epsilon_{n,p} \Big( \int\limits_{T_0}^{T_1} \|\partial_t R(t)\|_\infty^2 dt \Big) \int\limits_{-1}^1 |h_{n,p}|^2 dx 
\,,
\end{array}
\end{equation}
where $\|\partial_t R (t)\|_{\infty} := \|\partial_t R(t,.) \|_{L^\infty(-1,1)}$ and
$$|\epsilon_{n,p}| \leqslant \left\lbrace \begin{array}{l}
\epsilon_p'  := \Big( \frac{C_3}{(|p|+1)} + e^{C_3 \left( 2 + \frac{1}{T_1-T_0} - C_4 \min\{|p|\,,\,p^2\} (T_1-T_0-T_*) \right)} \Big),  \;\forall (n,p) \in\mathbb{Z}^2
\vspace{.2cm}\\
\epsilon_n'' := \Big( \frac{C_3}{n^2} + e^{C_3 \left( 1 + \frac{1}{T_1-T_0} - C_4 n^2 (T_1-T_0) \right)} \Big),\; \text{ if } |n| > 2|p|\,.
\end{array}\right.$$

\noindent \emph{Step 2: proof of the existence of a constant $C=C(T_1-T_0)>0$ such that, for $(n,p) \in \mathbb{Z}^2$ large enough, the following inequality holds}
$$  \frac{2}{R_0^2} \int\limits_{-1}^1 |\partial_t g_{n,p}(T_1,x)|^2 dx \leqslant 
 C \int\limits_{T_0}^{T_1} \int\limits_a^b |\partial_t g_{n,p}|^2 dx dt  +  \frac{1}{2} \int_{-1}^1 |h_{n,p}|^2 dx\,.$$ 
Note that $\epsilon_p' \longrightarrow 0$ when $|p| \rightarrow \infty$ and 
$\epsilon_n'' \longrightarrow 0$ when $|n| \rightarrow \infty$,
thus there exists $\rho>0$ such that
\begin{equation} \label{hyp_rho}
\frac{2 \max\{\epsilon_j'\,,\,\epsilon_j''\} }{R_0^2}   \int\limits_{T_0}^{T_1} \|\partial_t R(t)\|_\infty^2 dt < \frac{1}{2}\,, \quad\forall j \in \mathbb{Z} \text{ with } |j|>\rho\,.
\end{equation}
Let $(n,p) \in \mathbb{Z}^2$ be such that $n^2 + p^2 > 5 \rho^2$. 

\noindent \emph{\underline{First case}: $|p|>\rho$.} We have that
$$\frac{2 \epsilon_{n,p}}{R_0^2} \int\limits_{T_0}^{T_1} \|\partial_t R(t)\|_\infty^2 dt
\leqslant \frac{2 \epsilon_p'}{R_0^2} \int\limits_{T_0}^{T_1} \|\partial_t R(t)\|_\infty^2 dt < \frac{1}{2}\,.$$

\noindent \emph{\underline{Second case}: $|p| \leqslant \rho$.} Since $n^2 > 4 \rho^2$, we have  $|n|>2|p|$ and $|n|>\rho$. Then
$$\frac{2 \epsilon_{n,p}}{R_0^2} \int\limits_{T_0}^{T_1} \|\partial_t R(t)\|_\infty^2 dt
\leqslant \frac{2 \epsilon_n''}{R_0^2} \int\limits_{T_0}^{T_1} \|\partial_t R(t)\|_\infty^2 dt < \frac{1}{2}\,.$$
Step 2 follows with $C:=\frac{2}{R_0^2}\exp\left( C_3\left( 2 + \frac{1}{T_1-T_0} \right)\right)$
thanks to (\ref{interm}).

\noindent \emph{Step 3: conclusion.} For $(n,p) \in \mathbb{Z}^2$ such that $n^2 + p^2 > 5 \rho^2$,
we deduce from (\ref{interm1}) and Step 2 that
$$\frac{1}{2} \int\limits_{-1}^1 |h_{n,p}|^2 dx 
\leqslant C  \int\limits_{T_0}^{T_1} \int\limits_a^b |\partial_t g_{n,p}|^2 dx dt
+ \frac{2}{R_0^2} \int_{-1}^1 |\hl _{n,p} g_{n,p}(T_1,x)|^2 dx\,. \eqno{\Box}$$

\section{3D-Lipschitz stability estimate when $\omega$ is a tube}
\label{Sec:tube}

The goal of this section is the proof of Theorem \ref{Thm:Lips_Stab}.
\\

For $n, p \in \mathbb{Z}$, $H _{n,p}:=L^2(-1,1) \otimes e^{i(ny+pz)}$ is a closed subspace of $L^2(\Omega)$.
For $j \in \mathbb{N}$, we define
$$E_{j,p}:=\oplus_{|n| \leqslant 2^{j}} H_{n,p}$$
and denote by $\Pi_{j,p}$ the orthogonal projection from $L^2(\Omega)$ onto $E_{j,p}$.  
We also denote by $\Pi_{\infty,p}$ the orthogonal projection from $L^2(\Omega)$ onto $L^2((-1,1)\times\mathbb{T}) \otimes e^{ipz}$.
Moreover, $Id$ stands for the identity operator on $L^2(\Omega)$.

\subsection{Observability with source for frequency packets}

The goal of this section is the proof of the following result.

\begin{Prop} \label{Prop:packets}
Let $a,b \in \mathbb{R}$ be such that $-1 \leqslant a < b \leqslant 1$ and let $C_3, C_4>0$ be as in Proposition \ref{Prop:1D-IO-source}.
Let $\omega_y$ be an open subset of $\mathbb{T}$ and $\omega:=(a,b) \times \omega_y \times \mathbb{T}$.
There exists $C_5>C_3$ and $C_6\in (0,C_4)$ such that,
for every $T>0$, $p,j \in \mathbb{Z}$ with 
\begin{equation} \label{def0:j0}
j \geqslant j_0(p) :=
\begin{cases}
\left[ \frac{\ln|p|}{\ln(2)} \right]+2
&  \text{ if } \quad p \neq 0
\vspace{.1cm}
\\
0 & \text{ if } \quad p=0\,.
\end{cases}
\end{equation}
$g^0 \in L^2(\Omega)$, $\widetilde{h} \in L^2((0,T)\times\Omega)$, the solution of \eqref{H} satisfies
\begin{multline*}
\int_{\Omega} | \Pi_{j,p} g(T) |^2dxdydz    \leqslant 
 e^{C_5 \left( 2^j + \frac{1}{T} - C_6 |p| T \right)} \int_0^T \int_\omega | \Pi_{j,p} g |^2dxdydz 
 \\
 \displaystyle +  \left(C_3 + e^{C_3\left( 1 + \frac{1}{T} + |p| - C_4 |p| T \right)} \right) \int_0^T \int_\Omega | \Pi_{j,p} \widetilde{h} |^2dxdydz\,. 
\end{multline*}
\end{Prop}

The proof of this result relies on the following spectral inequality.
 
\begin{Prop} \label{Lem:LR}
Let $\omega_y$ be an open subset of $\mathbb{T}$.
There exists $C_{LR}>0$ such that, for all $N \in \mathbb{N}^*$ and $(b_k)_{-N \leqslant k \leqslant N} \in \mathbb{C}^{2N+1}$,
$$\Big( \sum_{k=-N}^N |b_k|^2 \Big)^{1\over2} \leqslant e^{C_{LR} N} \Big( \int_{\omega_y} \Big| \sum\limits_{k=-N}^N b_k e^{i k y} \Big|^2 dy \Big)^{1\over2}\,.$$
\end{Prop}

In this statement, the functions $y \mapsto e^{iky}/\sqrt{2\pi}$ are the orthonormal eigenfunctions of the
Laplace operator on the 1D-torus $\mathbb{T}$.
In arbitrary dimension, for a second-order symmetric elliptic operator, typically
the Laplace-Beltrami operator $\Delta_g$ on a bounded Riemannian manifold $\mathcal{M}$ of dimension
$d$, with or without boundary, the spectral inequality takes the form
\begin{equation} \label{spect_ineq}
\|u\|_{L^2(\mathcal{M})} \leqslant C e^{C \sqrt{\mu}} \|u\|_{L^2(\omega)}\,, \quad u \in \text{Span}\{ \phi_j ; \mu_j \leqslant \mu \}\,,
\end{equation}
where $\omega \subset \mathcal{M}$ is an open subset of $\mathcal{M}$
and the functions $\phi_j$ form a Hilbert basis of $L^2(\mathcal{M})$ of eigenfunctions of $-\Delta_g$,
associated with the non negative eigenvalues $\mu_j$, $j \in \mathbb{N}$,
counted with their multiplicities.
(In the case of a manifold with boundary, one can consider homogeneous Dirichlet or Neuman boundary conditions).
This was proven in \cite{Lebeau-Robbiano, Lebeau_Jerison, Lebeau_Zuazua}.
\\

Inequality (\ref{spect_ineq}) is a key tool to prove the null controllability 
of the heat equation by the Lebeau-Robbiano strategy (see \cite{Lebeau-LeRousseau} for a presentation). 
This strategy was adapted much later to the case of separated variables, for the null controllability of parabolic equations in
stratified media in \cite{BDLR}: in one direction, one has observability by means of
a Carleman estimate for a one-dimensional parabolic operator with parameter, and,
in the transverse direction, a spectral inequality such as (\ref{spect_ineq}) is used. This
approach was successfully transposed to the study of the null controllability of the
Grushin equation in \cite{Grushin} and the Kolmogorov equation in \cite{MR3163490}.
This approach was also adapted to the study of Lipschitz stability for the Grushin equation in \cite{MR3162108}.
The strategy we develop in this article
is more subtle than the one above.
Indeed, the choice of the space variables with respect to which we develop in Fourier series is not arbitrary.
For instance, the strategy would not work by developing only with respect to $z$
because the 2D resulting heat equations would not satisfy appropriate Carleman estimates.
This is why we take the Fourier series with respect to both $y$ and $z$.
Then, we apply the Lebeau-Robbiano strategy with respect to $(y,n)$, 
paying attention to the behaviour of the different constants with respect to $p$ (the Fourier frequency associated with $z$).
Indeed, these constants need to be uniform with respect to $p$ to get the conclusion.
\\

\noindent \textbf{Proof of Proposition \ref{Prop:packets}:} 
Let $p \in \mathbb{Z}$ and $j \geqslant j_0(p)$, i.e., $2^j \geqslant 2^{j_0} > 2|p|$.
By the Bessel-Parseval equality, (\ref{1D-IO-source}), (\ref{def:epsilon_np_n<2p}) and the previous Lemma, we get
\begin{eqnarray*}
\lefteqn{\int_{\Omega} | \Pi_{j,p} g(T) |^2dxdydz =  \sum\limits_{ |n| \leqslant 2^j} \int_{-1}^1 |g_{n,p}(T,x)|^2 dx dt }
\\
&\leqslant&
\sum\limits_{ |n| \leqslant 2^j}   \Big[ e^{C_3 \left( 1 + \frac{1}{T} + |p| - C_4 |p| T \right)}   \int\limits_0^T \int\limits_a^b |g_{n,p}(t,x)|^2 dx dt  
\\
&  & + \Big( \frac{C_3}{|p|+1} + e^{C_3\left(1+\frac{1}{T}+|p|-C_4|p|T \right)} \Big) 
                                                         \int\limits_0^T \int\limits_{-1}^1 |\widetilde{h}_{n,p}(t,x)|^2 dx dt \Big]
\\
& \leqslant &  e^{C_3 \left( 1 + \frac{1}{T} + |p| - C_4 |p| T \right) + C_{LR} 2^j} 
\int\limits_0^T \int\limits_a^b \int\limits_{\omega_y} \left| \sum\limits_{ |n| \leqslant 2^j }  g_{n,p}(t,x) e^{iny} \right|^2 dx dy  dt 
\\
& & 
+  \Big( C_3 +  e^{C_3 \left(1+\frac{1}{T}+|p|-C_4|p|T \right)} \Big) 
\sum\limits_{ |n| \leqslant 2^j} \int\limits_0^T \int\limits_{-1}^1 |\widetilde{h}_{n,p}(t,x)|^2 dx dt
\\
&\leqslant &  e^{C_5 \left( 2^{j} + \frac{1}{T} - C_6 |p| T \right)} \int\limits_0^T \int\limits_\omega | \Pi_{j,p} g |^2 dxdydzdt
\\
& &
+ \Big( C_3 + e^{C_3\left(1+\frac{1}{T}+|p|-C_4|p|T \right)} \Big) \int\limits_0^T \int\limits_\Omega | \Pi_{j,p} \widetilde{h} |^2dxdydzdt\,,
\end{eqnarray*}
where $C_5 := 2 C_3 + C_{LR}$ and $C_6:=\frac{C_3 C_4}{C_5}$. $\hfill \Box$

\subsection{Lebeau-Robbiano strategy for high frequencies}

The goal of this section is the proof of the following result.

\begin{Prop} \label{Prop:rec}
There exists $C_7>0$ 
such that, for all $T \geqslant 1$, $p \in \mathbb{Z}$, $g^0 \in L^2(\Omega)$, and $\widetilde{h} \in L^2((0,T)\times \Omega)$
the solution of \eqref{H} satisfies 
$$\| (Id-\Pi_{j_0,p})g(T) \|^2 \leqslant 
C_7 T^2 \Big(  \int\limits_0^T \int\limits_\omega |\Pi_{\infty,p}g|^2 dxdydzdt
+ \int\limits_0^T \| \Pi_{\infty,p} \widetilde{h}(t) \|^2 dt \Big)$$
where $j_0=j_0(p)$ is as in (\ref{def0:j0}).
\end{Prop}

\begin{rk}\rm
The lower bound $T \geqslant 1$ is chosen arbitrarily and may be replaced by any positive lower bound $T \geqslant T_*>0$ with a constant $C_7=C_7(T_*)>0$.
In the proof, assuming $T \geqslant 1$ will simplify the expression of the $T$-dependance of several constants.
Such an assumption is compatible with the fact that the positive result we have in mind only holds in large time.
\end{rk}

To prove Proposition \ref{Prop:rec}, we follow the Lebeau-Robbiano strategy, from the observability point of view,
with respect to parameter $n$ keeping parameter $p$ fixed.
We pay attention to the dependence of constants with respect to $p$.
\\

In the whole section, we fix $\rho \in (0,1)$, $T \geqslant 1$, $p \in \mathbb{Z}$ and $j_0:=j_0(p)$ as in (\ref{def0:j0}).
Note that
\begin{equation} \label{def:j0}
2^{j_0-1} \leqslant 2 |p| < 2^{j_0} \quad \text{ if } p \neq 0\,.
\end{equation}
Let $K=K(T,p,\rho)>0$ be such that
$$T=2K\sum\limits_{j \geqslant j_0} 2^{-\rho j} = \frac{2K 2^{-\rho j_0}}{1-2^{-\rho}}.$$
From (\ref{def:j0}) it follows that
\begin{equation} \label{encadrement_K}
\frac{2^\rho-1}{2} T |p|^\rho < K\leqslant 2^\rho \frac{2^\rho-1}{2} T |p|^\rho \quad \text{ if } p \neq 0\,.
\end{equation}
Then there exists $K_*=K_*(\rho)>0$, independent of $(T,p)$, such that
\begin{equation} \label{def:K*}
 K(T,p,\rho) \geqslant 2 K_* T > 0\,,  \quad\forall (T,p) \in (0,\infty) \times \mathbb{Z}\,.
\end{equation}
We now define times
\begin{equation} \label{def:tauj}
\tau_j=\tau_j(T,p,\rho)=K 2^{-j\rho} \quad \text{ and  } \quad \alpha_j=\alpha_j(T,p,\rho)= 2 \sum_{k=j_0}^j  \tau_k \quad \forall j \geqslant j_0\,,
\end{equation}
and time intervals
$$I_j:=(T-\alpha_{j-1}-\tau_j,T-\alpha_{j-1}) \quad \text{ and } \quad J_j:=(T-\alpha_{j},T-\alpha_{j-1})\, \quad \forall j > j_0\,.
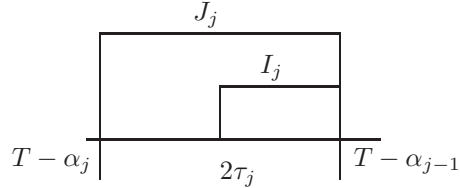
$$

\bigskip
\begin{figure}[h]
\begin{picture}(450,40)
\put(215,-10){\line(0,1){55}} 
\put(125,45){\line(1,0){90}}
\put(170,25){\line(1,0){45}}
\put(170,25){\line(0,-1){20}}
\put(92,-5){$T-\alpha_j$}
\put(220,-5){$T-\alpha_{j-1}$} 
\put(160,50){$J_j$}
\put(185,30){$I_j$}
\put(170,-10){$2\tau_j$}
\put(125,-10){\line(0,1){55}}
\thicklines 
\put(120,5){\line(1,0){110}}
\end{picture}
\end{figure}
\noindent
We will also use the notation 
\begin{equation} \label{def:lambda(2j)}
\lambda(2^j) = \frac{2^{2j}}{4} 
\end{equation}
so that $\lambda_{n,p} \geqslant \lambda(2^{j})$ for every $|n| \geqslant 2^{j}$ and $j \geqslant j_0(p)$ by (\ref{dissip_speed_n}) and (\ref{def:j0}).

We will need the following preliminary result, which is a consequence of  the Bessel-Parseval identity and Lemma \ref{Lem:Duhamel_np}.

\begin{Lem} \label{Lem:Duhamel_paquet}
Let $T_1, T_2 \in \mathbb{R}$,  $p \in \mathbb{Z}$, $j_1, j_2 \in \mathbb{N} \cup \{\infty\}$ be such that 
$$0 \leqslant T_1 < T_2 < \infty\quad\mbox{and}\quad j_0(p) \leqslant j_1 < j_2 \leqslant \infty.$$
For every $g^0 \in L^2(\Omega)$ and $\widetilde{h} \in L^2((0,T)\times\Omega)$, the solution of \eqref{H} satisfies
\begin{eqnarray*}
\| (\Pi_{j_2,p}-\Pi_{j_1,p})g(T_2) \|^2 \leqslant
2 \| (\Pi_{j_2,p}-\Pi_{j_1,p})g(T_1) \|^2 e^{- 2 \lambda(2^{j_1})(T_2-T_1)}
\\
+ \frac{1}{\lambda(2^{j_1})} \int_{T_1}^{T_2} \| (\Pi_{j_2,p}-\Pi_{j_1,p}) \widetilde{h}(t) \|^2 dt\,.
\end{eqnarray*}
\end{Lem}

\begin{Prop} \label{Prop:Prel_rec}
There exist $C_8, C_9 >0$ such that 
for every $T \geqslant 1$, $p \in \mathbb{Z}$, $j > j_0(p)$, $g^0 \in L^2(\Omega)$, and $\widetilde{h} \in L^2((0,T)\times \Omega)$
the solution of \eqref{H} satisfies
\begin{multline}
 e^{- C_8 2^j} \| \Pi_{j,p} g (T-\alpha_{j-1}) \|^2 \leqslant
 \int\limits_{I_j \times \omega} |\Pi_{\infty,p} g|^2 dxdydzdt
 \\
 + \frac{C_9 T}{2^j}  \int\limits_{J_j \times \Omega} | \Pi_{\infty,p} \widetilde{h}|^2dxdydzdt
 + T e^{- K_*T 2^{(2-\rho)j}} \| \Pi_{\infty,p} g(T-\alpha_j) \|^2,
\end{multline}
where $K_*=K_*(\rho)>0$ is as in \eqref{def:K*}.
\end{Prop}

\noindent \textbf{Proof of Proposition \ref{Prop:Prel_rec}:} 
Let $p \in \mathbb{Z}$, $j > j_0(p)$, $g^0 \in L^2(\Omega)$, $\widetilde{h} \in L^2((0,T)\times \Omega)$.
To simplify notations in this proof, we assume that $g_0 \in L^2((-1,1)\times\mathbb{T}) \otimes e^{ipz}$ and
$\widetilde{h} \in L^2(0,T;L^2((-1,1)\times\mathbb{T})\otimes e^{ipz})$, so that $\Pi_{\infty,p}g(t)=g(t)$ and $\Pi_{\infty,p} \widetilde{h}(t)=\widetilde{h}(t)$
for every $t \in [0,T]$. We also write $\Pi_j$ instead of $\Pi_{j,p}$ and omit all integration symbols such as $dx, dy, dz, dt$.

By Proposition \ref{Prop:packets}, the solution of \eqref{H} satisfies
\begin{multline}\label{obs_n_interm0}
 \| \Pi_{j} g(T-\alpha_{j-1}) \|^2 \leqslant
 e^{C_5 \left(2^j + \frac{1}{\tau_j} \right)} \int\limits_{I_j \times \omega} | \Pi_{j} g |^2 
 \\
 +\Big( C_3 + e^{C_3\left(1 + \frac{1}{\tau_j}+|p|-C_4|p|\tau_j \right)}\Big) \int\limits_{I_j \times \Omega} | \Pi_{j} \widetilde{h} |^2 \,.
\end{multline}
Moreover, we have
\begin{equation} \label{obs_n_interm1}
\int\limits_{I_j \times \omega} | \Pi_{j} g|^2 
 \leqslant 2 \int\limits_{I_j \times \omega} |  g |^2 +  2 \int\limits_{I_j \times \Omega} |(Id-\Pi_j) g|^2\,,\\
\end{equation}
and, by Lemma \ref{Lem:Duhamel_paquet} applied with $T_1=T-\alpha_j$, $T_2=t \in I_j$, $j_1=j$, $j_2=\infty$,
\begin{multline}\label{obs_n_interm2}
\int\limits_{I_j \times \Omega} |(Id-\Pi_j)g|^2 
\leqslant
2 \tau_j \|(Id-\Pi_j)g(T-\alpha_j)\|^2 e^{-2\lambda(2^j) \tau_j} 
 \\
 + \frac{\tau_j}{\lambda(2^j)} \int\limits_{J_j} \|(Id-\Pi_j) \widetilde{h} \|^2\,.
\end{multline}
Therefore,
\begin{multline}\label{Packet_interm}
  \| \Pi_{j} g(T-\alpha_{j-1}) \|^2
   \leqslant  2  e^{C_5 \left(2^j + \frac{1}{\tau_j} \right)} \int\limits_{I_j \times \omega} |g|^2
 \\
+ 2  e^{C_5 \left(2^j + \frac{1}{\tau_j} \right)} 
\Big( \frac{\tau_j}{\lambda(2^j)} 
+ C_3 e^{-C_5\left(2^j + \frac{1}{\tau_j} \right)} 
+ e^{C_3\left(1+\frac{1}{\tau_j}+|p|\right)-C_5\left(2^j+\frac{1}{\tau_j} \right)} \Big) 
\int\limits_{J_j} \| \widetilde{h} \|^2 
\\
+ 4 \tau_j e^{C_5 \left(2^j + \frac{1}{\tau_j} \right)-2\lambda(2^j) \tau_j} \|(Id-\Pi_j)g(T-\alpha_j)\|^2.
\end{multline}
From (\ref{def:tauj}), (\ref{def:K*}),  assumptions $T\geqslant 1$ and $\rho < 1$, we deduce that
$$\frac{1}{\tau_j} = \frac{2^{j\rho}}{K} \leqslant \frac{2^{j\rho}}{K_* T} \leqslant \frac{2^j}{K_*}\,, \quad \forall j>j_0(p)\,.$$
Then there exists $C_8>0$ independent of $(T,p,g^0,\widetilde{h})$ such that
\begin{equation} \label{step1:C1}
2  e^{C_5 \left(2^j + \frac{1}{\tau_j} \right)} \leqslant e^{C_8 2^j}\,, \quad \forall j>j_0(p)\,.
\end{equation}
We also have 
\begin{eqnarray*}
\lefteqn{C_3 \Big(1+\frac{1}{\tau_j}+|p|\Big)-C_5 \Big(2^j+\frac{1}{\tau_j} \Big)}
\\
&= & C_3 - \left( C_5 - C_3 \right) \frac{1}{\tau_j} - C_3 \left( 2^j - |p| \right) - \left( C_5 - C_3 \right) 2^j
\\
&\leqslant & C_3 -  \left( C_5 - C_3 \right) 2^j  \text{ because } C_5 > C_3 \text{ and } 2|p| < 2^{j_0} < 2^j\,.
\end{eqnarray*}
Thus, there exists a constant $c>0$ independent of $(T,p,g^0,\widetilde{h})$ such that
$$e^{C_3\left(1+\frac{1}{\tau_j}+|p|\right)-C_5\left(2^j+\frac{1}{\tau_j} \right)}
\leqslant e^{C_3 -  \left( C_5 - C_3 \right) 2^j } \leqslant \frac{c}{2^j}\,, \quad \forall j \geqslant j_0(p)\,,$$
$$\frac{\tau_j}{\lambda(2^j)} \leqslant \frac{2 T}{2^{2j}} \leqslant \frac{cT}{2^j}\,, \quad \forall j \geqslant j_0(p)\,,$$
and
$$C_3 e^{-C_5\left(2^j + \frac{1}{\tau_j} \right)} \leqslant C_3 e^{-C_5 2^j} \leqslant \frac{c}{2^j} \,, \quad \forall j \geqslant j_0(p)\,.$$
As a consequence, there exists $C_9>0$ independent of $(T,p,g^0,\widetilde{h})$ such that
\begin{equation} \label{step2:C2}
\frac{\tau_j}{\lambda(2^j)} 
+ C_3 e^{-C_5\left(2^j + \frac{1}{\tau_j} \right)} 
+ e^{C_3\left(1+\frac{1}{\tau_j}+|p|\right)-C_5\left(2^j+\frac{1}{\tau_j} \right)}
\leqslant \frac{C_9 T }{2^j}\,, \quad \forall j>j_0(p)\,.
\end{equation}
By (\ref{def:lambda(2j)}), (\ref{def:tauj}) and (\ref{def:K*}), we have
\begin{equation} \label{step3:C3} 
4 \tau_j e^{-2 \lambda(2^j) \tau_j}
\leqslant T e^{-K_* T 2^{(2-\rho)j}} 
\,, \quad \forall j>j_0(p)\,,
\end{equation}
because $4 \tau_j \leqslant 2 (\tau_j + \tau_{j_0}) \leqslant T$.
Finally, from (\ref{Packet_interm}), (\ref{step1:C1}), (\ref{step2:C2}) and (\ref{step3:C3}) we deduce that
\begin{multline*}
 \| \Pi_j  g(T-\alpha_{j-1}) \|^2 \leqslant
 e^{C_8 2^j} \int\limits_{I_j \times \omega} |g|^2
 + e^{C_8 2^j} \frac{C_9 T }{2^j}  \int_{J_j \times \Omega} | \widetilde{h} |^2 
 \\
 + T e^{C_8 2^j - K_* T 2^{(2-\rho)j}} \|(Id-\Pi_j) g(T-\alpha_j)\|^2,
\end{multline*}
which ends the proof of Proposition \ref{Prop:Prel_rec}. $\hfill \Box$
\\

\noindent \textbf{Proof of Proposition \ref{Prop:rec}:} 
Let $p \in \mathbb{Z}$, $j > j_0(p)$, $g^0 \in L^2(\Omega)$ and let $\widetilde{h} \in L^2((0,T)\times \Omega)$.
To simplify notations in this proof, we assume that $g_0 \in L^2((-1,1)\times\mathbb{T}) \otimes e^{ipz}$ and
$\widetilde{h} \in L^2(0,T;L^2((-1,1)\times\mathbb{T})\otimes e^{ipz})$, so that $\Pi_{\infty,p}g(t)=g(t)$ and $\Pi_{\infty,p} \widetilde{h}(t)=\widetilde{h}(t)$
for every $t \in [0,T]$. We also write $\Pi_j$ instead of $\Pi_{j,p}$  and omit all integration symbols such as $dx, dy, dz, dt$.
Let $C_8$, $C_9$ be as in Proposition \ref{Prop:Prel_rec}.
\\

\noindent \emph{Step 1: we prove by induction on $j \geqslant j_0+1$ that, for every $j \geqslant j_0+1$,}
\begin{multline*}
 \sum\limits_{k=j_0+1}^{j}  e^{-C_8 2^k} \| \Pi_k  g(T-\alpha_{k-1})\|^2
\hfill (\mathcal{P}_j)
 \\
  \leqslant 
 \sum\limits_{k=j_0+1}^{j} \delta_k \int\limits_{I_k \times \omega} |g|^2
+ A_j \int_{T-\alpha_j}^{T}  \| \widetilde{h}(t) \|^2 dt  + B_j \| g(T-\alpha_j) \|^2
\end{multline*}
{\em where}
\begin{equation} \label{ABe_n1}
\delta_{j_0+1}:=1, \quad
A_{j_0+1}:= \frac{C_9 T}{2^{j_0+1}}, \quad 
B_{j_0+1}:= T e^{-K_* T 2^{(2-\rho)(j_0+1)}}
\end{equation}
{\em and }
\begin{equation} \label{rec:delta}
\delta_{j+1}:=1+B_j e^{C_8 2^{j+1}},
\end{equation}
\begin{equation} \label{rec:A}
A_{j+1}:=  A _j +  \frac{B_j}{2^{2j}}  + \frac{\delta_{j+1} C_9 T}{2^{j+1}} ,
\end{equation}
\begin{equation} \label{rec:B}
B_{j+1}:= ( 2 B_{j} + \delta_{j+1}  T ) e^{-K_* T 2^{(2-\rho)(j+1)}}.
\end{equation}
The inequality $(\mathcal{P}_{j_0+1})$ is given by Proposition \ref{Prop:Prel_rec} with $j=j_0+1$. 
Let us now assume that $(\mathcal{P}_j)$ holds for some $j>j_0$ and prove $(\mathcal{P}_{j+1})$.
We have
$$B_j \| g(T-\alpha_j) \|^2 =  B_j \| \Pi_{j+1} g(T-\alpha_j)\|^2 +  B_j \|(Id-\Pi_{j+1}) g(T-\alpha_j)\|^2.$$
Applying Lemma \ref{Lem:Duhamel_paquet} to the last term
(with $T_1=T-\alpha_{j+1}$, $T_2=T-\alpha_{j}$, $j_1=j+1$, $j_2=\infty$)
we get 
\begin{multline*}
 B_j \| g(T-\alpha_j) \|^2  \leqslant B_j \| \Pi_{j+1} g(T-\alpha_j)\|^2 
 \\
 + 2 B_j \|(Id-\Pi_{j+1}) g(T-\alpha_{j+1})\|^2 e^{-4 \lambda(2^{j+1}) \tau_{j+1}}
 + \frac{B_j}{\lambda(2^{j+1})}  \int\limits_{J_{j+1}}  \| (Id-\Pi_{j+1}) \widetilde{h}(t) \|^2 dt \,.
\end{multline*}
Moreover, by (\ref{def:tauj}), (\ref{def:lambda(2j)}) and (\ref{def:K*}), we have 
$$4 \lambda(2^{j+1}) \tau_{j+1} = K 2^{(2-\rho)(j+1)}  \geqslant K_* T 2^{(2-\rho)(j+1)}
\quad \text{ and } \quad
\frac{B_j}{\lambda(2^{j+1})}=\frac{B_j}{2^{2j}}\,.$$
Thus, $(\mathcal{P}_j)$ implies
\begin{eqnarray}\label{rec:interm1}
\lefteqn{ \sum\limits_{k=j_0+1}^{j}  e^{-C_8 2^k} \| \Pi_k  g(T-\alpha_{k-1})\|^2 - B_{j} \| \Pi_{j+1} g(T-\alpha_j)\|^2}
 \\
 \nonumber
 & &
 \leqslant  \sum\limits_{k=1}^{j} \delta_k \int\limits_{I_k \times \omega} |g|^2   + \Big( A_j + \frac{B_j}{2^{2j}} \Big) \int_{T-\alpha_{j+1}}^T
 \| \widetilde{h} \|^2+ 2 B_j e^{-K_* T 2^{(2-\rho)(j+1)}} \|g(T-\alpha_{j+1})\|^2.
\end{eqnarray}
Moreover, by Proposition \ref{Prop:Prel_rec}, we also have
\begin{multline}\label{rec:interm2}
 e^{-C_8 2^{j+1}} \| \Pi_{j+1} g (T-\alpha_{j}) \|^2 \leqslant
  \int\limits_{I_{j+1} \times \omega} |g|^2   + \frac{C_9 T}{2^{j+1}} \int_{J_{j+1}} \| \widetilde{h} \|^2
 \\
+ T e^{-K_* T 2^{(2-\rho)(j+1)}} \| g(T-\alpha_{j+1}) \|^2.
\end{multline}
Note that $\delta_{j+1}$ is chosen so that
$$\delta_{j+1} e^{-C_8 2^{j+1}} - B_j = e^{-C_9 2^{j+1}}.$$
Thus, summing (\ref{rec:interm1}) and $\delta_{j+1}*$(\ref{rec:interm2}), we get $(\mathcal{P}_{j+1})$, which ends the first step.
\\

\noindent \emph{Step 2: existence of $B^*>0$ independent of $(T,p) \in [1,\infty)\times\mathbb{Z}$ such that}
$$\widetilde{B}_j:=B_j e^{C_8 2^{j+1}} \leqslant B^* T\,, \quad \forall j>j_0\,.$$
From (\ref{rec:B}), (\ref{rec:delta}) and assumption \textquotedblleft $T\geqslant 1$ \textquotedblright, we deduce that 
$$\widetilde{B}_{j+1} \leqslant 3 T \left( \widetilde{B}_{j} +  1 \right) e^{C_8 2^{j+2} - K_* T 2^{(2-\rho)(j+1)}}\,, \quad \forall j>j_0\,.$$
Moreover, there exists $M_1, M_2>0$ independent of $(T,p) \in [1,\infty)\times\mathbb{Z}$ such that
$$ 3 T e^{-\frac{K_*}{4} T 2^{(2-\rho)(j+1)}} \leqslant 3 T e^{-\frac{K_*}{4} T 2^{(2-\rho)}} \leqslant M_1\,, \quad \forall j>j_0\,,$$
and
$$ e^{C_8 2^{j+2} - \frac{K_*}{4} T 2^{(2-\rho)(j+1)}} \leqslant e^{C_8 2^{j+2} - \frac{K_*}{4} 2^{(2-\rho)(j+1)}} \leqslant M_2\,, \quad \forall j>j_0\,.$$
Thus 
\begin{equation} \label{rec_Bjtilde_bis}
\widetilde{B}_{j+1} \leqslant M \left(  \widetilde{B}_{j} + 1 \right) e^{-\frac{K_*}{2} T 2^{(2-\rho)(j+1)}}\,, \quad \forall j>j_0\,.
\end{equation}
where $M:=M_1 M_2$ is independent of $(T,p) \in [1,\infty)\times\mathbb{Z}$ and may be assumed to be $>1$. In particular
$$\widetilde{B}_{j+1} \leqslant M \left(  \widetilde{B}_{j} + 1 \right)\,, \quad \forall j>j_0\,,$$
or
$$\widetilde{B}_{j+1} + \frac{M}{M-1} \leqslant M \Big(  \widetilde{B}_{j} + \frac{M}{M-1} \Big)\,, \quad \forall j>j_0\,.$$
Thus,
$$\widetilde{B}_{j} + 1 \leqslant \widetilde{B}_{j} + \frac{M}{M-1}  \leqslant M^{j-j_0-1} \Big(  \widetilde{B}_{j_0+1} + \frac{M}{M-1} \Big)$$
and by (\ref{rec_Bjtilde_bis})  we deduce that
$$\widetilde{B}_{j+1} \leqslant  M^{j-j_0} \Big( \widetilde{B}_{j_0+1}  + \frac{M}{M-1} \Big) e^{-\frac{K_*}{2} T 2^{(2-\rho)(j+1)}}\,, \quad \forall j>j_0\,.$$
Moreover, there exists $c_1>0$, independent of $(T,p) \in [1,\infty)\times\mathbb{Z}$, such that
$$M^{j-j_0} e^{-\frac{K_*}{2} T 2^{(2-\rho)(j+1)}} \leqslant M^{j-j_0} e^{-\frac{K_*}{2} 2^{(2-\rho)(j+1)}} \leqslant c_1\,, \quad \forall j>j_0\,$$
and, in view of \eqref{def:j0},
$$\widetilde{B}_{j_0+1} = T e^{C_8 2^{j_0+2} - K_* T 2^{(2-\rho)(j_0+1)}} \leqslant T e^{16 C_8 |p|  - K_* |p|^{2-\rho} } \leqslant c_1 T\,.$$
This ends Step 2, because $T \geqslant 1$.
\\

\noindent \emph{Step 3: existence of $A^*>0$ independent of $(T,p) \in [1,\infty)\times\mathbb{Z}$ such that} 
$$A_j \leqslant A^* T^2\,, \quad \forall j>j_0\,.$$
By definition, we have
\begin{equation*}
A_{j}  = A_{j_0+1}+\sum\limits_{k=j_0+1}^{j-1} \Big( \frac{B_k}{2^{2k}}  + \frac{\delta_{k+1} C_9 T }{2^{k+1}} \Big)
\leqslant \frac{C_9 T }{2^{j_0+1}} + \sum_{k=0}^\infty \Big( \frac{B^* T }{2^{2k}} + \frac{\left(1+B^*T\right) C_9 T}{2^{k+1}} \Big)
\end{equation*}
which proves Step 3, because $T \geqslant 1$.
\\

\noindent \emph{Step 4: passing to the limit as $j \rightarrow \infty$ in $(\mathcal{P}_j)$.}
The last term on the right-hand side of $(\mathcal{P}_j)$ converges to zero because  $B_j \leqslant B^* T e^{-C_8 2^{j+1}}$.
Thus, we get 
\begin{multline}\label{sum_infinie}
  \sum\limits_{k=j_0+1}^{\infty}  e^{-C_8 2^k}  \| \Pi_k g(T-\alpha_{k-1}) \|^2
 \\
\leqslant \left(1+B^* T \right)   \int\limits_0^T \int\limits_{\omega} |g|^2 +  A^* T^2 \int\limits_0^T \|\widetilde{h}(t)\|^2 dt\,. 
\end{multline}
\noindent \emph{Step 5: conclusion.}
Using the Pythagorean theorem and Lemma \ref{Lem:Duhamel_paquet}, we get
\begin{eqnarray*}
\lefteqn{ \| (Id-\Pi_{j_0})g(T) \|^2 =    \sum_{k=j_0+1}^{\infty} \|(\Pi_k-\Pi_{k-1})g(T)\|^2}
 \\
& \leqslant &
 2 \sum_{k=j_0+1}^\infty   \|(\Pi_k-\Pi_{k-1})g(T-\alpha_{k-1})\|^2 e^{-2\lambda(2^{k-1})\alpha_{k-1}}
 \\ & & \hspace{2.5cm} + \sum_{k=j_0+1}^\infty \frac{1}{\lambda(2^{k-1})} \int_{T-\alpha_{k-1}}^{T} \|\widetilde{h}\|^2\,.
\end{eqnarray*}
Moreover, there exists $c_2>0$ independent of $(T,p)\in [1,\infty)\times\mathbb{Z}$ such that
$$e^{-2\lambda(2^{k-1})\alpha_{k-1}} \leqslant e^{- 2 \lambda(2^{k-1})\tau_{k-1}} \leqslant e^{- K_* 2^{(2-\rho)(k-1)}} \leqslant \frac{c_2}{2} e^{-C_8 2^k}\,, \forall k>j_0\,.$$
Thus
\begin{eqnarray*}
\lefteqn{\| (Id-\Pi_{j_0})g(T) \|^2}
\\
& \leqslant&
c_2 \sum\limits_{k=j_0+1}^\infty   \| \Pi_k g(T-\alpha_{k-1}) \|^2 e^{-C_8 2^k} 
+ \Big( \sum_{k=0}^\infty \frac{1}{2^{2k}} \Big) \int_{0}^{T} \|\widetilde{h}\|^2\,.
\end{eqnarray*}
Finally (\ref{sum_infinie}) gives the conclusion of Proposition \ref{Prop:rec} because $T \geqslant 1$. $\hfill \Box$

\subsection{3D observability inequality with source term}

The goal of this section if the proof of the following result.

\begin{Prop} \label{Prop:Obs_3D_ac_source}
There exist $T_*>0$ and  $c_*:(T_*,\infty) \rightarrow (0,\infty)$ continuous such that,
for every $T>T_*$, $p \in \mathbb{Z}$, $g^0 \in L^2(\Omega)$, and $\widetilde{h} \in L^2((0,T)\times \Omega)$, 
the solution of \eqref{H} satisfies 
$$\| \Pi_{\infty,p} g(T) \|^2 \leqslant c_*(T) \Big( \int\limits_0^T \int\limits_\omega | \Pi_{\infty,p} g|^2 
+ \int\limits_0^T \int\limits_\Omega | \Pi_{\infty, p} \widetilde{h}|^2  \Big)\,.$$
\end{Prop}

\noindent \textbf{Proof of Proposition \ref{Prop:Obs_3D_ac_source}:}
Let $p \in \mathbb{Z}$, $g^0 \in L^2(\Omega)$, and $\widetilde{h} \in L^2((0,T)\times \Omega)$.
To simplify notations in this proof, we assume that $g_0 \in L^2((-1,1)\times\mathbb{T}) \otimes e^{ipz}$ and
$\widetilde{h} \in L^2(0,T;L^2((-1,1)\times\mathbb{T})\otimes e^{ipz})$, so that $\Pi_{\infty,p}g(t)=g(t)$ and $\Pi_{\infty,p} \widetilde{h}(t)=\widetilde{h}(t)$
for every $t \in [0,T]$. We also write $\Pi_j$ instead of $\Pi_{j,p}$.
Let $C_4>0$ be as in Proposition \ref{Prop:1D-IO-source}, $C_6\in (0,C_4)$ be as in Proposition \ref{Prop:packets} and $T_*:=\max\{1\,,\,8/C_6\}$. We assume that $T>T_*$.
By orthogonality,
\begin{equation} \label{Pyth}
\|g(T)\|^2 = \|\Pi_{j_0} g(T)\|^2 + \|(Id-\Pi_{j_0})g(T)\|^2\,.
\end{equation}
Appealing to Proposition \ref{Prop:packets}, we get
\begin{multline*}
\int\limits_{\Omega} | \Pi_{j_0} g(T) |^2   \leqslant
e^{C_5 \left( 2^{j_0} + \frac{2}{T} - C_6|p| \frac{T}{2} \right)} \int\limits_{T/2}^T \int\limits_\omega | \Pi_{j_0} g |^2 
\\
 +  \Big( C_3 + e^{C_3\left( 1 + \frac{2}{T} + |p| - C_4 |p| \frac{T}{2} \right)} \Big) \int\limits_{T/2}^T \int\limits_\Omega | \Pi_{j_0} \widetilde{h} |^2\,. 
\end{multline*}
Moreover, invoking (\ref{def:j0}) and the fact that $T > T_* \geqslant \frac{8}{C_6} > \frac{8}{C_4}$, we obtain
$$2^{j_0}  - C_6 |p| \frac{T}{2} \leqslant 4 |p| - C_6 |p| \frac{T}{2} < 0 
\quad \text{ and } \quad 
|p| - C_4 |p| \frac{T}{2} < 0\,.$$
Thus, recalling  that $T \geqslant 1$ once again, we conclude that
\begin{equation} \label{Pck_simplified}
\int\limits_{\Omega} | \Pi_{j_0} g(T) |^2    \leqslant    
e^{2 C_5} \int\limits_{T/2}^T \int\limits_\omega | \Pi_{j_0} g |^2  
+  \left( C_3 + e^{3 C_3} \right) \int\limits_{T/2}^T \int\limits_\Omega | \Pi_{j_0} \widetilde{h} |^2\,.
\end{equation}
By Lemma \ref{Lem:Duhamel_paquet}, we have
\begin{multline}\label{Duh}
  \int\limits_{T/2}^T \int\limits_\omega | \Pi_{j_0} g |^2   
  \leqslant
2 \int_{T/2}^T \int_\omega | g |^2  + 2 \int_{T/2}^T \int_\Omega | (Id-\Pi_{j_0}) g |^2
 \\
\leqslant
2 \int_{T/2}^T \int_\omega | g |^2  + 2 T \Big\|(Id-\Pi_{j_0}) g\Big(\frac{T}{2}\Big) \Big\|^2 
+ T \int_{T/2}^T \|\widetilde{h}\|^2\,.
\end{multline}
Therefore,  (\ref{Pyth}), (\ref{Pck_simplified}), and (\ref{Duh}) yield
\begin{multline*}
\|g(T)\|^2 \leqslant  
2 e^{ 2 C_5 } \int_{T/2}^T \int_\omega | g |^2  
+ \left(  T  e^{ 2C_5 }  +   C_3  + e^{3 C_3} \right)  \int_{T/2}^T \|\widetilde{h}\|^2 
 \\
 +\|(Id-\Pi_{j_0}) g(T)\|^2 + 2 T e^{2C_5} \Big\|\left(Id-\Pi_{j_0}\right) g\Big(\frac{T}{2}\Big) \Big\|^2\,.
\end{multline*}
We complete the proof by applying Proposition \ref{Prop:rec} to the last two terms. $\hfill \Box$

\subsection{Proof of Theorem \ref{Thm:Lips_Stab}}

Let $T_*$ be as in Proposition \ref{Prop:Obs_3D_ac_source} and let $T_0 \in [0,T_1)$ be such that $T_1-T_0>T_*$.
In view of \eqref{Hyp:R}, we have
$$\int_\Omega |h|^2 \leqslant \frac{1}{R_0^2} \int_{\Omega} |R(T_1) h|^2 
\leqslant \frac{2}{R_0^2} \int_{\Omega} \left( |\partial_t g(T_1)|^2 + |\hl g(T_1)|^2 \right)\,.$$
By the Bessel-Parseval identity (note the particular form of $\omega=(a,b)\times\omega_y\times\mathbb{T}$)
and Proposition \ref{Prop:Obs_3D_ac_source}, we obtain
$$\int_\Omega |h|^2 
\leqslant \frac{2 C_{10}(T)}{R_0^2} \Big( \int_{T_0}^{T_1} \!\!\!\int_{\omega} |\partial_t g |^2 
+ \int_{T_0}^{T_1}  \|\partial_t R(t)\|_\infty^2 \|h\|^2 dt \Big)
+ \frac{2}{R_0^2} \int_\Omega |\hl g(T_1)|^2$$
for some constant $C_{10}>0$. The conclusion follows with
$$\eta(T):=\frac{R_0}{2 \sqrt{C_{10}(T)}}\,. \eqno{\Box}$$

\section{3D-Observability inequality when $\omega$ is a tube}
\label{Sec:Obs}

The goal of this section is the proof of Theorem \ref{Thm:Obs}.

\subsection{Observability in large time} \label{subsec:Observability in large time}

Let $T_*$ be as in Proposition \ref{Prop:Obs_3D_ac_source} and $T>T_*$.
The observability of \eqref{H} on $\omega=(a,b) \times \omega_y \times \mathbb{T}$ in time $T>T_*$ follows from the Bessel-Parseval identity and Proposition \ref{Prop:Obs_3D_ac_source}
(no source term $\widetilde{h}$).

\subsection{No observability in small time}

The goal of this section is the proof of the following result.

\begin{Prop} \label{Prop:Obs_negatif}
Let $a,b \in \mathbb{R}$ be such that $-1 < a < b < 1$ and $$\omega:=(a,b) \times \mathbb{T} \times \mathbb{T}.$$
If $T<\frac{1}{8} \max\{ (1+a)^2 , (1-b)^2 \}$, then  \eqref{H} is not observable in $\omega$ in time $T$.
\end{Prop}

\noindent \textbf{Proof of Proposition \ref{Prop:Obs_negatif}:} One may assume that $-1 < a < b=1$. 
Let $T<\frac{1}{8} (1+a)^2$. We are going to construct a sequence $(g_k)_{k \in \mathbb{N}^*}$ of solutions  of \eqref{H}
such that
\begin{equation} \label{c-ex}
\dfrac{\int_0^T \int_a^1 \int_{\mathbb{T}} \int_\mathbb{T} |g_k(t,x,y,z)|^2 dz dy dx dt}{\int_{-1}^1 \int_\mathbb{T} \int_\mathbb{T} |g_k(T,x,y,z)|^2 dz dy dx dt}
\underset{k \rightarrow \infty}{\longrightarrow} 0\,.
\end{equation}
Let $\alpha:=\frac{1-a}{2}>0$ and $\epsilon>0$ be such that 
\begin{equation}\label{hyp:epsilon}
(- 1 + \alpha)^2 - \epsilon >0\,, \quad 
T < \frac{1}{8}(a+1)^2 - \epsilon 
\end{equation}
and $k_1(\epsilon) \in \mathbb{N}^*$ be such that
\begin{equation} \label{hyp:k*}
\left( \pm 1+\frac{[\alpha k]}{k} \right)^2 \geqslant \left( \pm 1+\alpha \right)^2 - \epsilon\,, \forall k \geqslant k_1(\epsilon)\,,
\end{equation}

\noindent
\emph{Step 1: construction of $g_k$ from an explicit approximate solution.}
The function 
$$G(x):=\frac{1}{\sqrt[4]{\pi}} e^{-\frac{x^2}{2}}$$
satisfies
$$\left\lbrace \begin{array}{l}
-G''(x)+x^2 G(x)=G(x)\,, \quad  x \in \mathbb{R}\,, \\
\int_{\mathbb{R}} G(x)^2 dx =1\,.
\end{array}\right.$$
Let $\theta_{\pm} \in C^\infty_c(\mathbb{R})$ be such that 
$$\theta_{\pm}(\pm 1)=1\,, \quad  \theta_{\pm}(\mp 1)=0 \quad  \text{ and } \quad  \text{Supp}(\theta_-) \cap (a,1) = \emptyset\,.$$
For $(n,p) \in \mathbb{Z} \times \mathbb{R}^*_+$, the function
$$\mathcal{K}_{n,p}(t,x):=\sqrt[4]{p} \left\{  G\left(\sqrt{p}\Big(x+\frac{n}{p}\Big)\right) 
- \sum\limits_{\sigma \in \{-1,1\}} G\left(\sqrt{p}\Big(\sigma+\frac{n}{p}\Big)\right) \theta_{\sigma}(x)  \right\} e^{-pt}$$
satisfies
$$\left\lbrace \begin{array}{ll}
\Big( \partial_t - \partial_x^2 + (px+n)^2 \Big) \mathcal{K}_{n,p}(t,x)=E_{n,p}(t,x)\,, \quad & (t,x) \in (0,\infty)\times(-1,1)\,, \\
\mathcal{K}_{n,p}(t,\pm 1)=0\,,                                                               & t \in (0,\infty)\,,
\end{array}\right.$$
where
$$E_{n,p}(t,x)=\sqrt[4]{p} \sum\limits_{\sigma \in \{-1,1\}} \Big(-p+\partial_x^2 -(px+n)^2 \Big) \theta_{\sigma}(x) e^{-pt}  G\left(\sqrt{p}\Big(\sigma+\frac{n}{p}\Big)\right)\,.$$
For $(n,p) \in \mathbb{Z} \times \mathbb{R}^*_+$, let $\mathcal{G}_{n,p}(t,x)$ be the solution of
$$\left\lbrace \begin{array}{ll}
\Big( \partial_t - \partial_x^2 + (px+n)^2 \Big) \mathcal{G}_{n,p}(t,x)=0\,, \quad & (t,x) \in (0,\infty)\times(-1,1)\,, \\
\mathcal{G}_{n,p}(t,\pm 1)=0\,,                                                    & t \in (0,\infty)\,,\\
\mathcal{G}_{n,p}(0,x)=\mathcal{K}_{n,p}(0,x)\,,                                   & x \in (-1,1)\,.
\end{array}\right.$$
Then, by Duhamel's formula, there exists $c_1>0$, independent of $(n,p) \in \mathbb{Z}\times\mathbb{R}^*_+$, such that
for all $(t,n,p) \in (0,T)\times\mathbb{Z}\times\mathbb{N}^*$
$$\|(\mathcal{G}_{n,p}-\mathcal{K}_{n,p})(t,.)\|_{L^2(-1,1)}^2 \leqslant c_1 \int_0^t \|E_{n,p}(s)\|_{L^2(-1,1)}^2 ds\,.$$
Thus, recalling the definition of $E_{n,p}$ we conclude that there exists $c_2>0$, independent of $(n,p) \in \mathbb{Z}\times\mathbb{R}^*_+$, such that
for every $(t,n,p) \in (0,T)\times\mathbb{Z}\times\mathbb{R}^*_+$
\begin{equation} \label{erreur}
\|(\mathcal{G}_{n,p}-\mathcal{K}_{n,p})(t,.)\|_{L^2(-1,1)} \leqslant c_2 \frac{p^2+n^2}{\sqrt[4]{p}} 
\max_{\sigma \in \{-1,1\}}  e^{-\frac{p}{2}\left(\sigma+\frac{n}{p}\right)^2} \,.
\end{equation}
We define
$$g_k(t,x,y,z):= \mathcal{G}_{[\alpha k],k}(t,x) e^{i([\alpha k] y+ k z)}\,.$$

\noindent
\emph{Step 2: we start the proof of \eqref{c-ex} arguing by contradiction.}
Assume that there exists $c_3>0$ such that, for every $k \in \mathbb{N}^*$,
$$\Big(\int\limits_{-1}^1\! \int\limits_\mathbb{T} \!\int\limits_\mathbb{T} |g_k(T,x,y,z)|^2 dz dy dx \Big)^{1\over2}
\leqslant c_3 \Big( \int\limits_0^T \!\int\limits_a^1\! \int\limits_\mathbb{T}\! \int\limits_{\mathbb{T}} |g_k(t,x,y,z)|^2 dz dy dx dt \Big)^{1\over2}\,.$$
Thanks to the Bessel-Parseval identity, the above inequality may be written as
$$\Big( \int\limits_{-1}^1  \left| \mathcal{G}_{[\alpha k],k}(T,x) \right|^2  dx \Big)^{1\over2}
\leqslant c_3 \Big( \int\limits_0^T \!\int\limits_a^1 \left| \mathcal{G}_{[\alpha k],k}(t,x) \right|^2  dx dt \Big)^{1\over2}
\quad \forall k \in \mathbb{N}^*\,.$$
By the triangular inequality and (\ref{erreur}), we deduce that, for some constant $c_4>0$,  
\begin{multline}\label{c-ex-1}
 \Big( \int\limits_{-1}^1  |\mathcal{K}_{[\alpha k],k}(T,x)|^2  dx \Big)^{1\over2}
\leqslant  
  c_3 \Big( \int\limits_0^T \int\limits_a^1 \left| \mathcal{K}_{[\alpha k],k}(t,x) \right|^2  dx dt \Big)^{1\over2} \\
+ c_4 k^{7/8}  e^{-\frac{k}{2}\left[(-1+\alpha)^2 - \epsilon \right]}\,,
\quad \forall k > k_1(\epsilon).
\end{multline}

\noindent
\emph{Step 3: lower bound for the left-hand side of \eqref{c-ex-1}.}
We have
\begin{multline*}
 \Big( \int\limits_{-1}^1  \Big| \mathcal{K}_{[\alpha k],k}(T,x) \Big|^2  dx \Big)^{1\over2} \geqslant 
 \Big( \int\limits_{-1}^1 \sqrt{k} e^{-k\big(x+\frac{[\alpha k]}{k}\big)^2} e^{-2 k T}   dx \Big)^{1\over2} 
 \\
 - \sum\limits_{\sigma \in \{-1,1\}} \Big( \int\limits_{-1}^1  
\sqrt{k} e^{- k \big( \sigma+\frac{[\alpha k]}{k} \big)^2} \theta_{\sigma}(x)^2 e^{-2 k T}   dx \Big)^{1\over2}\,.
\end{multline*}
Thus, there exists $c_5, c_6, c_7>0$ and $k_2(\epsilon) \geqslant k_1(\epsilon)$ such that, for all $ k>k_2(\epsilon)$,
\begin{multline*}
 \Big( \int\limits_{-1}^1  \left| \mathcal{K}_{[\alpha k],k}(T,x) \right|^2  dx \Big)^{1\over2}
 \geqslant 2 c_5 e^{-k T}   
            - c_6   e^{-\frac{k}{2}\left[(-1+\alpha)^2 - \epsilon + 2T \right]}  
 \\
  \geqslant e^{-kT } \left( 2 c_5 - c_6 e^{-\frac{k}{2}\left[(-1+\alpha)^2-\epsilon\right]}  \right) 
   \geqslant c_7 e^{-kT}, 
  \end{multline*}
where we have also taken \eqref{hyp:epsilon} into account.

\smallskip
\noindent
\emph{Step 4: upper bound for the right-hand side of \eqref{c-ex-1}.}
There exist constants $c_8, c_9, c_{10}>0$ and $k_3(\epsilon)>k_2(\epsilon)$ such that, for every $k>k_3(\epsilon)$,
\begin{eqnarray*}
\lefteqn{\Big( \int\limits_0^T \!\int\limits_a^1  \left| \mathcal{K}_{[\alpha k],k}(t,x) \right|^2   dx dt \Big)^{1\over2}}
\\
& \leqslant & 
\Big( \int\limits_0^T\! \int\limits_a^1 \sqrt{k} e^{-k\left(x+\frac{[\alpha k]}{k}\right)^2} e^{-2 k t}  dx dt \Big)^{1\over2}
 +
\Big( \int\limits_0^T \!\int\limits_a^1 \sqrt{k} e^{-k\left(1+\frac{[\alpha k]}{k}\right)^2} \theta_{+}(x)^2 e^{-2 k t} dx dt \Big)^{1\over2}
\\
& \leqslant &
\Big(   \int\limits_{\sqrt{k}\left(a+\frac{[\alpha k]}{k}\right)}^\infty  \frac{e^{-x^2}}{2k} dx   \Big)^{1\over2}
+ \frac{c_8}{\sqrt[4]{k}} e^{-\frac{k}{2}\left[(1+\alpha)^2 - \epsilon\right]}
\\ &\leqslant &
c_9 \left(   \frac{e^{-k\left(a+\frac{[\alpha k]}{k}\right)^2}}{\sqrt{k}\left(a+\frac{[\alpha k]}{k}\right)}   \right)^{1\over2}
+ c_8 e^{-\frac{k}{2}\left[(1+\alpha)^2 - \epsilon\right]} \leqslant 
c_{10} e^{-\frac{k}{2}\left[(a+\alpha)^2 - \epsilon\right]}
\end{eqnarray*}
where we have used the fact that $0 < a + \alpha < 1+\alpha$.

\smallskip
\noindent
\emph{Step 5: conclusion.} We deduce from (\ref{c-ex-1}), Step 3 and Step 4 that
$$c_7 e^{- kT } \leqslant c_3 c_{10} e^{-\frac{k}{2}\left[(a+\alpha)^2 - \epsilon\right]}
+ c_4 k^{7/8} e^{-\frac{k}{2}\left[(-1+\alpha)^2 - \epsilon \right] }\,, \;\;\forall k>k_3(\epsilon)\,.$$
Moreover, by choice of $\alpha$, we have $(a+\alpha)^2=(-1+\alpha)^2$, thus
$$e^{-kT} \leqslant [\ c_3 c_{10} + c_4 k^{7/8}] e^{-\frac{k}{2}\left[(-1+\alpha)^2 - \epsilon\right]}\,, \forall k>k_3(\epsilon)\,.$$
This is a contradiction because $T<(-1+\alpha)^2 - \epsilon$. $\hfill \Box$

\subsection{Proof of Theorem \ref{Thm:Obs}}

Let $a,b \in \mathbb{R}$ be such that $-1 < a < b < 1$, $\omega_y$ be an open subset of $\mathbb{T}$ and $\omega:=(a,b) \times \omega_y \times \mathbb{T}$.
The quantity
$$T_{\min}:=\inf\left\{~T>0~:~\text{ system } \eqref{H} \text{ is observable in } \omega \text{ in time } T~\right\}$$
is finite by Section \ref{subsec:Observability in large time} and $\geqslant \frac{1}{8} \max\{ (1+a)^2 , (1-b)^2 \}$ by Proposition \ref{Prop:Obs_negatif}.

\section{Observability on an unbounded domain}
\label{sec:z_dans_R}

In this section, we consider the Heisenberg equation
\begin{equation} \label{H0}
\left\lbrace \begin{array}{ll}
\Big(\partial_t  - \left(\partial_{x_1}-\frac{x_2}{2} \partial_{x_3} \right)^2 - \left(\partial_{x_2}+\frac{x_1}{2} \partial_{x_3} \right)^2 \Big) G = 0 
& \mbox{in } (0,T) \times \widetilde{\Omega}\,, \\
G(t,\pm 1,x_2,x_3)=0\,,                                                                   &  \\
G(t,x_1,-\pi,x_3)=G(t,x_1,\pi,x_3)\,,        &\\                                             
\partial_y G(t,x_1,-\pi,x_3) = \partial_y G(t,x_1,\pi,x_3)\,,                        &     \\
G(0,x)=G_0(x)\,,                                                                          & 
\end{array} \right.
\end{equation}
where $(x_1,x_2,x_3)\in\widetilde{\Omega}:=(-1,1) \times (-\pi,\pi) \times \mathbb{R}$,
and we prove the following observability result.

\begin{thm} \label{Thm:Main0}
Let $-1<a<b<1$, $-\pi<c<d<\pi$ and $\omega:=(a,b) \times (c,d) \times \mathbb{R}$.
Then there exists $T_{\min} \geqslant \frac{1}{8} \max\{ (1+a)^2 , (1-b)^2 \}$ such that
\begin{itemize}
\item for every $T>T_{\min}$, system \eqref{H0} is observable in $\omega$ in time $T$,
\item for every $T<T_{\min}$, system \eqref{H0} is not observable in $\omega$ in time $T$.
\end{itemize}
\end{thm}
\noindent
In a similar way, one can extend the Lipschitz stability result of Theorems~\ref{Thm:Lips_Stab_simple} and \ref{Thm:Lips_Stab}  to system \eqref{H0}.
\\

Observe that the change of variables
\begin{equation} \label{CVAR}
G(t,x_1,x_2,x_3)=g\left(t,x=x_1,y=x_2,z=x_3+\frac{x_1 x_2}{2}\right)
\end{equation}
transforms system (\ref{H0}) into the following auxiliary system
\begin{equation} \label{H_R}
\left\lbrace \begin{array}{ll}
\Big(\partial_t  - \partial_x^2 - (x \partial_z+\partial_y)^2 \Big) g(t,x,y,z) = 0\,, & (t,x,y,z) \in (0,T) \times \Omega\,, \\
g(t,\pm 1,y,z)=0\,,                                                                   & (t,y,z) \in (0,T)\times \mathbb{T} \times \mathbb{R}\,, \\
g(0,x,,z)=g_0(x,y,z)\,,                                                               & (x,y,z) \in \Omega\,.
\end{array} \right.
\end{equation}
where $\Omega=(-1,1)\times \mathbb{T}\times\mathbb{R}$.
This equation is well posed in $L^2(\Omega)$ as is equation (\ref{H0}).
Theorem \ref{Thm:Main0} is a direct consequence of the same statement for (\ref{H_R}).
The observability in large time can be proved by following the same arguments than in the previous sections,
replacing  summations over $p \in \mathbb{Z}$ by integrals over $p \in \mathbb{R}$.
In Section \ref{Sec:tube}, the assumption \textquotedblleft$p \in \mathbb{Z}$\textquotedblright \,  was used to simplify the writing of several estimates,
but the same analysis can be performed for $p \in \mathbb{R}$ by replacing $|p|$ by $\min\{|p|;p^2\}$ at several places,
as in Proposition \ref{Prop:1D-IO-source}.
On the other hand, the counter-example we gave to show that observability fails in time $T< \frac{1}{8} \max\{ (1+a)^2 , (1-b)^2 \}$
needs adjustment, which is what we do below.
\\

\noindent \textbf{Adaptation of the proof of Proposition \ref{Prop:Obs_negatif}:}
Let $k_1(\epsilon) \in \mathbb{N}^*$ be such that
\begin{equation} 
\Big( -1+\frac{k}{p} \Big)^2 \geqslant \left( -1+\alpha \right)^2 - \epsilon\,, \quad \forall p \in \Big( \frac{k}{\alpha}, 1+\frac{k}{\alpha} \Big)\,, \;\;\forall k \geqslant k_1(\epsilon).
\end{equation}

\noindent
\emph{Step 1: construction of $g_k$.} After introducing $G$, $\mathcal{K}_{n,p}$ and $\mathcal{G}_{n,p}$ as in the proof of Proposition \ref{Prop:Obs_negatif}, 
we define
$$g_k(t,x,y,z):=\frac{e^{iky}}{2\pi} \int_{\frac{k}{\alpha}}^{1+\frac{k}{\alpha}} \mathcal{G}_{k,p}(t,x) e^{ipz} dp \,.$$

\noindent
\emph{Step 2: contradiction argument.}
Suppose $c_1>0$ is such that, $\forall\,k \in \mathbb{N}^*$,
$$\Big(\int\limits_{-1}^1\! \int\limits_{\mathbb{T}}\! \int\limits_\mathbb{R} |g_k(T,x,y,z)|^2 dz dy dx \Big)^{1\over2}
\leqslant c_1\Big( \int\limits_0^T \!\int\limits_a^1 \!\int\limits_{\mathbb{T}}\! \int\limits_{\mathbb{R}} |g_k(t,x,y,z)|^2 dz dy dx dt \Big)^{1\over2}.$$
By Plancherel's identity, this inequality may be rewritten as
$$\Big( \int\limits_{-1}^1\! \int\limits_{\frac{k}{\alpha}}^{1+\frac{k}{\alpha}} |\mathcal{G}_{k,p}(T,x)|^2  dp dx \Big)^{1\over2}
\leqslant c_1 \Big( \int\limits_0^T \!\int\limits_a^1 \!\int\limits_{\frac{k}{\alpha}}^{1+\frac{k}{\alpha}} |\mathcal{G}_{k,p}(t,x)|^2  dp dx dt \Big)^{1\over2}
\quad \forall k \in \mathbb{N}^*\,.$$
As above, by the triangular inequality and (\ref{erreur}) we deduce that, for some constant $c_2>0$ and all $k > k_1(\epsilon)$,
\begin{multline}\label{c-ex-1_R}
  \Big( \int\limits_{-1}^1 \int\limits_{\frac{k}{\alpha}}^{1+\frac{k}{\alpha}} |\mathcal{K}_{k,p}(T,x)|^2  dp dx \Big)^{1\over2}
 \\
\leqslant
c_2 \Big( \int\limits_0^T \!\int\limits_a^1\! \int\limits_{\frac{k}{\alpha}}^{1+\frac{k}{\alpha}} |\mathcal{K}_{k,p}(t,x)|^2  dp dx dt \Big)^{1\over2}
+ c_2 k^{7/8}  e^{-\frac{k}{2\alpha}\left[(-1+\alpha)^2 - \epsilon \right]}.
\end{multline}

\noindent
\emph{Step 3: lower bound for the left-hand side of \eqref{c-ex-1_R}.}
We have
\begin{multline*}
\Big( \int\limits_{-1}^1\! \int\limits_{\frac{k}{\alpha}}^{1+\frac{k}{\alpha}} |\mathcal{K}_{k,p}(T,x)|^2  dp dx \Big)^{1\over2} \geqslant
 \Big( \int\limits_{-1}^1\! \int\limits_{\frac{k}{\alpha}}^{1+\frac{k}{\alpha}} \sqrt{p} e^{-p\left(x+\frac{k}{p}\right)^2} e^{-2pT}  dp dx \Big)^{1\over2} 
 \\
 - \sum\limits_{\sigma \in \{-1,1\}} \Big( \int\limits_{-1}^1\! \int\limits_{\frac{k}{\alpha}}^{1+\frac{k}{\alpha}} 
\sqrt{p} e^{-p\left(\sigma+\frac{k}{p}\right)^2} \theta_{\sigma}(x)^2 e^{-2pT}  dp dx \Big)^{1\over2}.
\end{multline*}
Thus, there exists $c_3, c_4>0$ and $k_2(\epsilon) \geqslant k_1(\epsilon)$ such that for all $k>k_2(\epsilon)$
\begin{eqnarray*}
\lefteqn{ \Big( \int\limits_{-1}^1 \int\limits_{\frac{k}{\alpha}}^{1+\frac{k}{\alpha}} |\mathcal{K}_{k,p}(T,x)|^2  dp dx \Big)^{1\over2}}
\\
& \geqslant & 2 c_3 \Big( \int\limits_{\frac{k}{\alpha}}^{1+\frac{k}{\alpha}} e^{-2pT} dp \Big)^{1/\over2}   
            - c_4  \Big( \int\limits_{\frac{k}{\alpha}}^{1+\frac{k}{\alpha}} 
                e^{-p\left[(-1+\alpha)^2 - \epsilon + 2T \right]} dp   \Big)^{1\over2} 
\\
& \geqslant & e^{-\frac{kT}{\alpha}} \left\{ 2 c_3 - c_4 \left( e^{-\frac{k}{2\alpha}[(-1+\alpha)^2-\epsilon]} \right) \right\}
                 \; \geqslant \;c_3\, e^{-\frac{kT}{\alpha}} ,               
\end{eqnarray*}
where we have also used \eqref{hyp:epsilon}.

\smallskip
\noindent
\emph{Step 4: upper bound for the right-hand side of \eqref{c-ex-1_R}.}
There exist constants $c_5, c_6, c_7>0$ and $k_3(\epsilon)>k_2(\epsilon)$ such that, for every $k>k_3(\epsilon)$,
\begin{eqnarray*}
\lefteqn{\Big( \int\limits_0^T \!\int\limits_a^1 \!\int\limits_{\frac{k}{\alpha}}^{1+\frac{k}{\alpha}} |\mathcal{K}_{k,p}|^2  dp dx dt \Big)^{1\over2}
\leqslant \Big( \int\limits_0^T \!\int\limits_a^1\! \int\limits_{\frac{k}{\alpha}}^{1+\frac{k}{\alpha}}
\sqrt{p} e^{-p\left(x+\frac{k}{p}\right)^2} e^{-2pt} dp dx dt \Big)^{1\over2}}
\\
&  & 
\hspace{3cm} +
\Big( \int\limits_0^T\! \int\limits_a^1\! \int\limits_{\frac{k}{\alpha}}^{1+\frac{k}{\alpha}}
\sqrt{p} e^{-p\left(1+\frac{k}{p}\right)^2} \theta_{+}^2 e^{-2pt} dp dx dt \Big)^{1\over2}
\\
& \leqslant &
\Big(  \int\limits_{\frac{k}{\alpha}}^{1+\frac{k}{\alpha}}\! \int\limits_{\sqrt{p}\left(a+\frac{k}{p}\right)}^\infty  \frac{e^{-x^2}}{2p} dx dp  \Big)^{1\over2}
+ \frac{c_5}{\sqrt[4]{k}} e^{-\frac{k}{2\alpha}\left[(1+\alpha)^2 - \epsilon\right]}
\\
& \leqslant &
c_6 \Big(  \int\limits_{\frac{k}{\alpha}}^{1+\frac{k}{\alpha}} \frac{e^{-p\left(a+\frac{k}{p}\right)^2}}{\sqrt{p}\left(a+\frac{k}{p}\right)} dp  \Big)^{1\over2}
+ c_5 e^{-\frac{k}{2\alpha}\left[(1+\alpha)^2 - \epsilon\right]}
\leqslant c_7 e^{-\frac{k}{2\alpha}\left[(a+\alpha)^2 - \epsilon\right]},
\end{eqnarray*}
where we have used the fact that $0< a + \alpha < 1+\alpha$.

\smallskip
\noindent
\emph{Step 4: conclusion.} Combining \eqref{c-ex-1_R}, Step 3, and Step 4 we conclude that
$$c_3 e^{-\frac{kT}{\alpha}} \leqslant c_2c_7 e^{-\frac{k}{2\alpha}\left[(a+\alpha)^2 - \epsilon\right]}
+ c_2 k^{7/8} e^{-\frac{k}{2\alpha}\left[(-1+\alpha)^2 - \epsilon \right] }$$
for all $k>k_3(\epsilon)$.
Moreover, the choice of $\alpha$ yields $(a+\alpha)^2=(-1+\alpha)^2$.
Thus, the above inequality gives a contradiction because $T<(a+\alpha)^2 - \epsilon$. $\hfill \Box$

\section{Conclusion and open problems} \label{Sec:ccl}

In this article, we have proved observability inequalities and Lipschitz stability estimates 
for the Heisenberg heat equation on product-shaped domains in $\R^3$.
Observations were taken on appropriate slices or tubes.
Both results require a  minimal time $T_{\min}>0$, a lower bound for which was given in terms of the distance between the observability region and the boundary of the space domain, in the $x$ direction.
The sharp evaluation of $T_{\min}>0$ is an open problem for which the techniques developed in \cite{minimal} for  Grushin's operator seem hard to utilize.

\smallskip
The Heisennberg heat equation is also well posed on the unbounded domain $(x,y,z) \in \mathbb{R} \times \mathbb{T} \times \mathbb{R}$.
In this situation, the dissipation spead $\lambda_{n,p}$ does no depend on $n$, because of the invariance under translations of variable $x$
(see Remark \ref{rk:validite_LR}). Thus the Lebeau-Robbiano method cannot be performed.
The validity of the observability inequality in this configuration is a completely open problem.

\appendix

\section{Proof of unique continuation}
\label{Appendix:UC}

In this appendix, we give a  proof of Proposition~\ref{Prop:Holmgren}. 
\\

Let $T>0$, $a,b\in \mathbb{R}$ be such that $-1<a<b<1$, $\omega_y$ be an open subset of $\T$ 
and $g \in C^0([0,T],L^2(\Omega))$ be a solution of \eqref{H} with $\tilde{h}=0$, which vanishes on $(0,T)\times(a,b)\times\omega_y\times\T$.
\\

Let $\epsilon>0$ be such that
\begin{equation}\label{eq:t_lowerbound}
\tilde{T}:=\frac{T}{\epsilon} > 2\,\mbox{diam}\,\big((-1,1)\times\T\times\T\big)
\end{equation}
and $\tilde{g}(\tau,x,y,z):=g(\epsilon \tau,x,y,z)$ for every $(\tau,x,y,z) \in (0,\tilde{T})\times\Omega$. Then
$$\left\lbrace \begin{array}{ll}
\left(\partial_{\tau}  - \epsilon \Big( \partial_x^2 + (x \partial_z+\partial_y)^2 \Big) \right) \tilde{g}(\tau,x,y,z) = 0\,, & (\tau,x,y,z) \in (0,\tilde{T}) \times \Omega\,, \\
\tilde{g}(\tau,\pm 1,y,z)=0\,,                                                                                      & (\tau,y,z) \in (0,\tilde{T})\times \mathbb{T} \times \mathbb{T}\,, \\
\tilde{g}(0,x,y,z)=g^0(x,y,z)\,,                                                                                    & (x,y,z) \in \Omega
\end{array} \right.$$
and $\tilde{g}=0$ on $(0,\tilde{T})\times(a,b)\times\omega_y\times\T$.
\\

Let $\mathcal{O}$ be the maximal open subset of $(0,\tilde{T}) \times (-1,1) \times \T \times \T$ such that $\tilde{g}=0$ on $\mathcal{O}$. Then
\begin{equation} \label{O_inclusion}
(0,\tilde{T}) \times (a,b) \times \omega_y \times \T\, \subset\, \mathcal{O}\,.
\end{equation}
Working by contradiction, we assume that $\mathcal{O} \neq (0,\tilde{T})\times(-1,1)\times \T \times \T$. Let 
\begin{equation} \label{def:txyz0}
\tau_0:=\frac{\tilde{T}}{2}\,, \quad  x_0:=\frac{a+b}{2}\,, \quad  y_0 \in \omega_y \quad \text{ and } \quad z_0 \in \T\,.
\end{equation}
Then, by (\ref{O_inclusion}), $(\tau_0,x_0,y_0,z_0)\in \mathcal{O}$. 
Let $(\tau_*,x_*,y_*,z_*) \in \partial \mathcal{O}$ be such that
$$\| (\tau_0,x_0,y_0,z_0) - (\tau_*,x_*,y_*,z_*) \| = r := \text{dist}\Big( (\tau_0,x_0,y_0,z_0) , \partial \mathcal{O} \Big)\,.$$
Then, necessarily
\begin{equation} \label{z0=z*}
z_0=z_*\,.
\end{equation}

\noindent
\emph{Step 1: we show that $\tau_* \in (0,\tilde{T})$.}
Working by contradiction, suppose that $\tau_* \in \{0,\tilde{T}\}$. 
Then, from (\ref{def:txyz0}) we deduce that $|\tau_0-\tau_*|=\tilde{T}/2$. So,
\begin{multline*}
 r = \| (\tau_0,x_0,y_0,z_0) - (\tau_*,x_*,y_*,z_*) \|
  \geqslant \tilde{T}/2 >  \mbox{diam}\,\big( (-1,1)\times\T\times\T \big)
 \\
   > \mbox{diam}\,\big( \mathcal{O} \cap [\{t_0\}\times(-1,1)\times\T\times\T ] \big)
    > \text{dist}\big( (\tau_0,x_0,y_0) , \partial \mathcal{O} \big)= r\,,
\end{multline*}
which is impossible.

\smallskip
\noindent
\emph{Step 2: we prove that}
\begin{equation} \label{normale_favorable}
\left(\begin{array}{c}
x_*-x_0  
\\
y_*-y_0
\end{array}\right)\neq 0.
\end{equation}
From (\ref{def:txyz0}), (\ref{O_inclusion}), and Step 1 we deduce that $(\tau_*,x_0,y_0,z_*)$ belongs to the open subset $\mathcal{O}$. 
So, $(\tau_*,x_0,y_0,z_*)\neq (\tau_*,x_*,y_*,z_*)$ since the latter point belongs to the boundary of $\mathcal{O}$.
Thus (\ref{normale_favorable}) holds.

\smallskip
\noindent
\emph{Step 3: we apply Holgren's uniqueness theorem.} 
We denote by 
$$\sigma((\tau,x,y,z),(s,\xi,\eta,\nu))=\epsilon\xi^2 + \epsilon(x\nu+\eta)^2$$ 
the principal symbol of the Heisenberg operator
$P:=\partial_{\tau} - \epsilon [ \partial_x^2 + (x \partial_z + \partial_y)^2]$.
Let $\Sigma$ be the sphere with center $(\tau_0,x_0,y_0,z_0)$ and radius $r$. 
By (\ref{z0=z*}), the unit normal to $\Sigma$ at $(\tau_*,x_*,y_*,z_*)$ is
\begin{equation*}
\vec{n}:=(n_{\tau},n_x,n_y,n_z)=\dfrac 1r \big(\tau_*-\tau_0,x_*-x_0,y_*-y_0,0\big).
\end{equation*}
Consequently,  $\sigma\big((\tau^*,x^*,y^*,z^*),(n_{\tau},n_x,n_y,n_z)\big)=n_x^2 + n_y^2  \neq 0$ by Step 2.
Thus $\Sigma$ is a smooth  noncharacteristic surface for $P$ at $(\tau_*,x_*,y_*,z_*)$.
Moreover, $g \equiv 0$ on one side of $\Sigma$, in a neighborhood of $(\tau_*,x_*,y_*,z_*)$.
By Holmgren's theorem~\cite[Theorem 8.6.5]{hormander},  $g \equiv 0$ on a open neighborhood of $(\tau^*,x^*,y^*,z^*)$. 
This contradicts the maximality of $\mathcal{O}$. \hfill $\Box$

\section{Carleman estimates for the 1D heat equation with parameters}

Let us set $\R_+=(0,\infty)$ and $\I=[-1,1]$. For any $T>0$ let $\I_T=[0,T]\times [-1,1]$.

\medskip
\noindent \textbf{Proof of Proposition \ref{Carleman_global}:}
Fix  $a',b'$ be such that $a<a'<b'<b$.
Fix a  real-valued function $\beta \in C^3([-1,1])$  such that
\begin{equation} \label{hyp-beta}
\beta \geqslant 1 \text{  on  } [-1,1])\,,
\end{equation}
\begin{equation} \label{hyp-beta'-(a,b)}
|\beta'|>0 \text{  on  } [-1,a'] \cup [b',1]\,,
\end{equation}
\begin{equation} \label{hyp_beta'-bord}
\beta'(1)>0\,,\quad \beta'(-1)<0\,,
\end{equation}
\begin{equation} \label{hyp-beta''}
\beta'' <0 \text{ on  }  [-1,a'] \cup [b',1]
\end{equation}
For any $M>0$ define
\begin{equation} \label{def-alpha}
\alpha(t,x)=\frac{M \beta (x)}{t(T-t)}\, ,\quad (t,x) \in (0,T) \times [-1,1]\,.
\end{equation}
Given a complex-valued function $g \in C^0([0,T];L^2(-1,1)) \cap L^2(0,T;H^1_0(-1,1))$, let us consider the standard transform
\begin{equation} \label{def-z}
z(t,x):=g(t,x) e^{-\alpha(t,x)}\, ,\quad (t,x) \in (0,T) \times [-1,1]\,.
\end{equation}
In the following computations we shall assume $g$ more regular so that we can compute derivatives of all the orders we need in order to obtain estimate \eqref{Carl_est}. Such a procedure can be made rigorous assuming $\mathcal{P}_{n,p} g\in L^2(\I_T)$.
We have
\begin{equation} \label{P123}
e^{-\alpha} \mathcal{P}_{n,p} g = P_{1}z + P_{2} z\, ,
\end{equation}
where we have set
\begin{equation} \label{def:P1-P2-P3}
\begin{array}{l}\displaystyle
P_{1}z= - \partial_x^2z +(\alpha_{t}-\alpha_{x}^{2}- \alpha_{xx}) z + (px+n)^2  z
\vspace{.2cm}\\
\displaystyle
P_{2}z= \partial_tz - 2 \alpha_{x} \partial_xz\, .
\end{array}
\end{equation}
We follow the classical proof 
which consists in taking the $L^{2}$-norm of both sides of the identity (\ref{P123}).
Developing the double product and recalling that $z$ is complex-valued, we obtain
\begin{equation} \label{P1P2<P3}
\int_{\I_T}{\mathcal Re}\big( P_{1}z\, \overline{P_{2}z} \big)\,  dx dt \leqslant \frac{1}{2} \int_{\I_T} | e^{-\alpha} \mathcal{P}_n g  |^2 dx dt\, ,
\end{equation}
where ${\mathcal Re}\,z $ denotes the real part of $z$.  We have
\begin{eqnarray*}
\lefteqn{{\mathcal Re}\big( P_{1}z\, \overline{P_{2}z} \big)=-{\mathcal Re}\big( \partial_x^2z\,\partial_t \overline{z} - 2 \alpha_{x}\,\partial_x^2z\, \partial_x \overline{z}\big)}
\\
&+ &
(\alpha_{t}-\alpha_{x}^{2}- \alpha_{xx})\,{\mathcal Re}\big( z\, \partial_t \overline{z} - 2 \alpha_{x}\,z\, \partial_x \overline{z}\big) +(px+n)^2{\mathcal Re}\big( z\,\partial_t \overline{z} - 2 \alpha_{x} \,z\,\partial_x \overline{z}\big)
\\
&=: &
Q_1+Q_2+Q_3\,.
\end{eqnarray*}
Now, we compute the integrals of $Q_1,Q_2$, and $Q_3$.

\smallskip
\noindent
\textbf{Evaluation of  $\int_{\I_T}Q_1$: }
integrating by parts, we get
\begin{eqnarray}\label{Q1}
\lefteqn{\int_{\I_T}Q_1\,dxdt}
\\
\nonumber
&=&\int_0^T\left[
 \alpha_{x} \left|\partial_xz\right|^2
 -{\mathcal Re}\left(\partial_xz\, \partial_t \overline{z}\right)\right]_{x=-1}^{x=1}dt
 \\
 \nonumber
 & &
 \hspace{4.5cm}
 +\int_{\I_T}\left[{\mathcal Re}\left(\partial_xz\,\partial_t\partial_x\overline{z}\right)
 - \alpha_{xx} \left|
 \partial_xz
 \right|^2\right]dxdt
 \\
\nonumber
&=& 
\int_{0}^{T} \left[
\alpha_{x}(t,1) \left|\partial_xz(t,1) \right|^{2} -
\alpha_{x}(t,-1)    \left|\partial_xz(t,-1)    \right|^{2} \right] dt 
-\int_{\I_T}
\alpha_{xx} 
\left|
\partial_xz
\right|^2 dxdt
\end{eqnarray}
because $\partial_t z(t,\pm 1)=0$ and $z(0,\cdot) \equiv z(T,\cdot) \equiv 0$.

\smallskip
\noindent
\textbf{Evaluation of  $\int_{\I_T}Q_3$: }
since $z(0,\cdot) \equiv z(T,\cdot) \equiv 0$ and $z(\cdot,-1) \equiv z(\cdot,1) \equiv 0$, we have
\begin{eqnarray}\label{Q3}
\lefteqn{\int_{\I_T}Q_3\,dxdt}
\\
\nonumber
&=&\int_{\I_T}(px+n)^2\left(\frac12 \partial_t|z|^2 - \alpha_{x} \partial_x|z|^2\right)dxdt
\\
\nonumber
&=&
\int_{\I_T}\big[(px+n)^2\alpha_{x}\big]_x |z|^2\,dxdt\,.
\end{eqnarray}

\smallskip
\noindent
\textbf{Evaluation of  $\int_{\I_T}Q_2$:} again integrating by parts, we have
\begin{eqnarray}\label{Q2}
\lefteqn{\int_{\I_T}Q_2\,dxdt}
\\
\nonumber
&=&\frac12\int_{\I_T}
(\alpha_{t}-\alpha_{x}^{2}- \alpha_{xx}) \partial_t|z|^2dxdt
-\int_{\I_T}
\alpha_x(\alpha_{t}-\alpha_{x}^{2}- \alpha_{xx}) \partial_x|z|^2dxdt
\\
\nonumber
&=&\int_{\I_T}
\left\{\left[\alpha_x(\alpha_{t}-\alpha_{x}^{2})\right]_x-\frac12\left(\alpha_{t}-\alpha_{x}^{2}- \alpha_{xx}\right)_t
- \left(\alpha_{x}^3\right)_x
\right\}|z|^2dxdt\,.
\end{eqnarray}

By combining  (\ref{Q1}), (\ref{Q3}),  and (\ref{Q2}) we obtain
\begin{eqnarray}\label{Q123}
\lefteqn{\int_{\I_T}{\mathcal Re}\left( P_{1}z\, \overline{P_{2}z} \right)dxdt}
 \\
 \nonumber
 &= &
-\int_{\I_T}
\alpha_{xx} \left(
\left|
\partial_xz
\right|^2 +3\alpha_{x}^2
|z|^2\right)dxdt
\\
\nonumber
& &+\int_{\I_T}
\left\{\left[\alpha_x(\alpha_{t}-\alpha_{xx})\right]_x-\frac12\left(\alpha_{t}-\alpha_{x}^{2}- \alpha_{xx}\right)_t
\right\}|z|^2dxdt
\\
\nonumber
& &+
\int_{\I_T}\left[(px+n)^2\alpha_{x}\right]_x |z|^2\,dxdt\,.
\end{eqnarray}
Now, observe that,  in view of \eqref{hyp-beta'-(a,b)} and \eqref{hyp-beta''},  
\begin{equation*}
m_1:=\min_{x\in [-1,a'] \cup [b',1]} |\beta'(x)|>0
\quad\mbox{and}\quad
m_2:=\min_{x\in [-1,a'] \cup [b',1]} -\beta''(x)>0
\end{equation*}
to deduce that
\begin{equation}\label{eq:bb}
-\alpha_{xx} \left(
\left|
\partial_xz
\right|^2 +3\alpha_{x}^2
|z|^2\right)\geqslant
\dfrac{m_2M}{t(T-t)} \left|\partial_xz\right|^2
+\dfrac{3m_1^2m_2M^3}{[t(T-t)]^3}|z|^2 
\end{equation}
for all $ x\in [-1,a'] \cup [b',1]$ and $t\in(0,T)$. Next, consider the function
\begin{equation}\label{eq:ra}
R_\alpha=\big[\alpha_x(\alpha_{t}-\alpha_{xx})\big]_x-\frac12\big(\alpha_{t}-\alpha_{x}^{2}- \alpha_{xx}\big)_t
\end{equation}
which is defined on $(0,T) \times [-1,1]$. Recalling \eqref{def-alpha}, one can easily check that
\begin{equation}\label{eq:|ra|}
|R_\alpha(t,x)|\leqslant \dfrac{C_0M^2}{[t(T-t)]^3}\,\|\beta\|_{\mathcal C^3(\I)}^2\,(T+T^2)\quad\forall (t,x)\in (0,T) \times [-1,1]
\end{equation}
for some constant $C_0>0$. Indeed, each of the terms that appear in \eqref{eq:ra} can be bounded by $M^2/[t(T-t)]^3$ times a polynomial of degree two with no zero order term in $\beta$ and its derivatives up to the third order, times $T$ or $T^2$.
Now,  for every
\begin{equation}\label{M1}
M \geqslant M_1(T,\beta):=\dfrac{C_0\|\beta\|_{\mathcal C^3(\I)}^2}{2m_1^2m_2}\,\,(T+T^2)\,,
\end{equation}
\eqref{eq:|ra|} implies that
\begin{equation*}
\Big(\dfrac{3m_1^2m_2M^3}{[t(T-t)]^3}+R_\alpha\Big)|z|^2 \geqslant \dfrac{m_1^2m_2M^3}{[t(T-t)]^3}|z|^2
\end{equation*}
for all $ x\in [-1,a'] \cup [b',1]$ and $t\in(0,T)$.
Therefore, owing to \eqref{Q123} and \eqref{eq:bb},  
\begin{multline} \label{In-Carl-2}
\displaystyle\int_{0}^{T}\!\!\! \int\limits_{(-1,a')\cup(b',1)}
\frac{C_{1} M}{t(T-t)} \left| \partial_xz \right|^{2} dxdt
\\
+ \displaystyle\int_{0}^{T}\!\!\! \int\limits_{(-1,a')\cup(b',1)}
\left\{
\frac{C_{3} M^{3}}{[t(T-t)]^3} |z|^{2}
+ [(px+n)^2 \alpha_x]_x |z|^2
\right\}
dxdt
\\
\displaystyle\leqslant 
\int_{0}^{T}\!\!\! \int_{a'}^{b'}
\left\{
\frac{C_{2} M}{t(T-t)}  \left| \partial_xz \right|^{2}
+ \frac{C_{4} M^{3}}{[t(T-t)]^3} |z|^{2}
- [(px+n)^2 \alpha_x]_x |z|^2
\right\} dxdt
\\
\displaystyle +  \int_{\I_T} | e^{-\alpha} \mathcal{P}_{n,p} g  |^2 dxdt
\end{multline}
 for some contants 
$C_j=C_j(\beta)>0\;(j=1,\dots,4)$.

Next, observe that, for every $x \in [-1,1]$
\begin{eqnarray}
\nonumber
 \big| \big[(px+n)^2 \alpha_x\big]_x \big| 
&=& \frac{M}{t(T-t)} \big| 2p(px+n)\beta'(x) + (px+n)^2 \beta''(x) \big|
\\
\label{In-Carl-3}
&\leqslant &\frac{C_5 M(n^2+p^2)}{t(T-t)}
\end{eqnarray}
where $C_5=C_5(\beta)>0$. Let
$$M_2=M_2(T,\beta,n,p):=\sqrt{\frac{2C_5}{C_3}} \left( \frac{T}{2} \right)^2 (|n|+p)$$
so that, for every $M \geqslant M_2$, we have
$$\frac{C_5 M(n^2+p^2)}{t(T-t)} \leqslant \frac{C_3 M^3}{2[t(T-t)]^3}\,.$$
From now on, we fix
$$M=\max\big\{ M_1(T,\beta)\,,\;M_2(T,\beta,n,p) \big\}$$
noting that, in view of \eqref{M1}, $M$ can be represented as in \eqref{M_Carleman} for some constant $\mathcal C_2(\beta)>0$. 
From \eqref{In-Carl-2} and \eqref{In-Carl-3}, it follows that
\begin{eqnarray}
\lefteqn{\int_{0}^{T} \!\!\!\int\limits_{(-1,a')\cup(b',1)} \left\{\frac{C_{1} M}{t(T-t)} \left| \partial_xz\right|^{2}  +  \frac{C_{3} M^{3}}{2[t(T-t)]^3} |z|^{2}
\right\}
}
\label{preCarl_z}
\\
\nonumber
&\leqslant &
\int_{0}^{T}\!\!\! \int_{a'}^{b'}
\left\{
\frac{C_{2} M}{t(T-t)}  \left| \partial_xz \right|^{2}
+ \frac{C_{6} M^{3}}{[t(T-t)]^3} |z|^{2}
\right\} dxdt
+  \int_{\I_T} | e^{-\alpha} \mathcal{P}_{n,p} g |^2 dxdt
\end{eqnarray}
where $C_6=C_6(\beta):=C_4+C_3/2$\,. 

A this point,  we need to recast the above inequality in terms of the original function $g$. Since, for every $\epsilon>0$,
\begin{multline*}
\frac{C_{1} M}{t(T-t)} \left| \partial_xg - \alpha_x g \right|^{2}
+
\frac{C_{3} M^{3}}{2 [t(T-t)]^{3}} |g|^{2}
\\
\geqslant
\Big( 1-\frac{1}{1+\epsilon} \Big) \frac{C_{1} M}{t(T-t)} \left| \partial_xg  \right|^{2}
+ \frac{M^3}{[t(T-t)]^{3}}\Big( \frac{C_{3}}{2} - \epsilon C_1 (\beta')^2 \Big) |g|^2\,,
\end{multline*}
choosing
$$\epsilon=\epsilon(\beta):=\frac{C_3}{4C_1\|\beta'\|_\infty^2}$$ 
from \eqref{preCarl_z} we deduce that
\begin{multline*} 
\int_{0}^{T}\!\!\!  \int\limits_{(-1,a') \cup (b',1)}  
\Big\{
\frac{C_{7} M}{t(T-t)} \left| \partial_xg \right|^{2}
+ \frac{C_{3} M^{3}|g|^{2}}{4[t(T-t)]^{3}} 
\Big\}\,  e^{-2\alpha} dxdt
\\
\leqslant 
\int_{0}^{T}\!\!\!  \int_{a'}^{b'} 
\Big\{
\frac{C_{9} M^{3} |g|^{2}}{[t(T-t)]^{3}}
+  
\frac{C_{8} M}{t(T-t)}  
\left| \partial_xg\right|^{2} 
\Big\}\,e^{-2\alpha} dxdt
+
 \int_{\I_T} | e^{-\alpha} \mathcal{P}_n g  |^2 dx dt
\end{multline*}
where 
\begin{eqnarray*}
C_7&=&C_7(\beta)=[1-1/(1+\epsilon)]C_1
\\
C_{8}&=&C_{8}(\beta)=2C_{2}
\\
C_{9}&=&C_{9}(\beta)=C_{6}+2C_{2} \sup \{ \beta'(x)^{2}:x \in [a',b'] \}.
\end{eqnarray*}
So, adding the same quantity to both sides, we obtain
\begin{eqnarray}
\lefteqn{\int_{\I_T}
\Big\{
\frac{C_{7} M}{t(T-t)} \left| \partial_xg \right|^{2}
+ \frac{C_{3} M^{3}|g|^{2}}{4[t(T-t)]^{3}} 
\Big\}\,  e^{-2\alpha} dxdt  }
\label{In-Carl-4}
\\
\nonumber
& &
\leqslant\int_{0}^{T}\!\!\! \int_{a'}^{b'} 
\Big\{
\frac{C_{11} M^{3} |g|^{2}}{(t(T-t))^{3}}
+  
\frac{C_{10} M}{t(T-t)}  
\left| \partial_xg\right|^{2} 
\Big\}\,e^{-2\alpha} dxdt
+ \int_{\I_T} | e^{-\alpha} \mathcal{P}_n g  |^2 dx dt
\end{eqnarray}
where $C_{10}=C_{10}(\beta)=C_{8}+C_7$ and $C_{11}=C_{11}(\beta)=C_{9}+C_3/4$.

The last step of the proof consists in showing that  $| \partial_xg|^{2}$ in the right-hand side of the above inequality can be absorbed  by the remaining two terms. This fact is a rather standard consequence of a Caccioppoli-type inequality. We give the proof for completenss. Let $\rho \in C^{\infty}(\mathbb{R})$ be such that $0\le \rho \le 1$ and
\begin{equation} \label{hyp-rho-hors(a,b)-carre}
\rho \equiv 1 \text{  on  } [a',b']\, ,
\end{equation}
\begin{equation} \label{hyp-rho-(a,b)-carre}
\rho \equiv 0 \text{  on  } [-1,a] \cup [b,1]\, .
\end{equation}
We have
$$\int_{\I_T} (\mathcal{P}_n g) \frac{g \rho e^{-2\alpha}}{t(T-t)} dx dt
= \int_0^T\!\!\! \int_{-1}^1 \left\{  
\partial_tg -\partial_x^2g + (px+n)^2 g
\right\}
\frac{g \rho e^{-2\alpha}}{t(T-t)} dx dt.$$
Integrating by parts with respect to time and space, we obtain
$$\int_{\I_T} \frac{1}{2} \partial_t g^2
\frac{\rho e^{-2\alpha}}{t(T-t)} dx dt
=\int_{\I_T} \frac{1}{2} |g|^2 \rho \Big\{
\frac{2\alpha_t}{t(T-t)} + \frac{T-2t}{[t(T-t)]^2}
\Big\}\, e^{-2\alpha} dx dt$$
and
\begin{multline*}
-\int_{\I_T}\partial_x^2g \frac{g \rho e^{-2\alpha}}{t(T-t)} dx dt
\\ 
= \int_{\I_T}
\frac{\rho e^{-2\alpha}}{t(T-t)}\left|\partial_xg\right|^2 dx dt- \int_{\I_T}  \frac{|g|^2 e^{-2\alpha}}{2t(T-t)} 
\big\{ \rho'' - 4 \rho'\alpha_x+\rho (4\alpha_x^2-2\alpha_{xx}) \big\}
 dx dt\,.
\end{multline*}
Thus,
\begin{multline}
\int_{\I_T} \mathcal{P}_n g \frac{g \rho e^{-2\alpha}}{t(T-t)} dx dt
\geqslant
\int_{\I_T} \frac{\rho e^{-2\alpha}}{t(T-t)} \left|\partial_xg\right|^2  dx dt
\\
- \int_{\I_T}  \frac{|g|^2 e^{-2\alpha}}{2t(T-t)} 
\Big\{ \rho'' - 4 \rho'\alpha_x+\rho \Big[ 4\alpha_x^2-2\alpha_{xx} -2\alpha_t-\frac{T-2t}{t(T-t)} \Big] \Big\} dx dt\,.
\end{multline}
Therefore,
\begin{eqnarray*}
\lefteqn{\int_0^T\!\!\! \int_{a'}^{b'} \frac{C_{10} M}{t(T-t)} \left| \partial_xg\right|^{2} e^{-2\alpha} dxdt
\leqslant
\int_{\I_T} \frac{C_{10} M \rho}{t(T-t)} \left| \partial_xg\right|^{2} e^{-2\alpha} dxdt}
\\
& \leqslant&
\int_{\I_T}
 \mathcal{P}_n g \frac{C_{10} M g \rho e^{-2\alpha}}{t(T-t)} dxdt
\\
&+ & \int_{\I_T}
\frac{C_{10} M |g|^2 e^{-2\alpha}}{2t(T-t)} 
\Big\{ \rho'' - 4 \rho'\alpha_x+\rho \Big[ 4\alpha_x^2-2\alpha_{xx} -2\alpha_t-\frac{T-2t}{t(T-t)} \Big] \Big\}
 dx dt
\\ &\leqslant&
\int_{\I_T} |\mathcal{P}_n g|^2 e^{-2\alpha} dx dt +
\int_0^T \!\!\!\int_a^b \frac{C_{12} M^3 |g|^2 e^{-2\alpha}}{[t(T-t)]^3} dx dt
\end{eqnarray*}
for some constant $C_{12}=C_{12}(\beta,\rho)>0$. Combining (\ref{In-Carl-4}) with the previous inequality, we get
\begin{multline*} 
\int_{\I_T} 
\Big\{
\frac{C_{7} M}{t(T-t)} \left| \partial_xg \right|^{2}
+ \frac{C_{3} M^{3}|g|^{2}}{4[t(T-t)]^{3}} 
\Big\}\,  e^{-2\alpha} dxdt
\\
\leqslant 
\frac{3}{2} \int_{\I_T}  | e^{-\alpha} \mathcal{P}_n g  |^2 dx dt
+
\int_{0}^{T}\!\!\! \int_{a}^{b} \frac{C_{13} M^{3} |g|^{2}}{[t(T-t)]^{3}} e^{-2\alpha} dxdt\, ,
\end{multline*}
where $C_{13}=C_{13}(\beta,\rho):=C_{11}+C_{12}$. 
Then, taking  
$$\mathcal{C}_1=\mathcal{C}_1(\beta):=\frac{\min\{C_7  ; C_3 /4 \}}{\max\{ 3/2 ; C_{13}  \}}$$
we obtain the global Carleman estimate (\ref{Carl_est}).
 \hfill $\Box$


\bibliography{biblio_heisenberg_1}
\bibliographystyle{plain}
\end{document}